\documentstyle[12pt]{article}

\newcommand{\sect}[1]{\setcounter{equation}{0}\section{#1}}

\newcommand{\app}[1]{\setcounter{section}{0}
\setcounter{equation}{0} \renewcommand{\thesection}{\Alph{section}}
\section{#1}}

\newcommand{\eq}{\begin{equation}}
\newcommand{\eqa}{\begin{eqnarray}}
\newcommand{\en}{\end{equation}}
\newcommand{\ena}{\end{eqnarray}}
\newcommand{\enn}{\nonumber \end{equation}}

\def\sk{\vskip .4cm}
\def\noi{\noindent}
\def\om{\omega}

\def\al{\alpha}
\def\la{\lambda}

\def\Ga{\Gamma}
\def\del{\delta}

\def\rinv{{1\over {r-r^{-1}}}}

\def\Cb{\bar{C}}

\def\fb{\bar{f}}

\def\rhop{{\rho}^{\prime}}

\def\thetap{{\theta}^{\prime}}

\def\onehalf{{1 \over 2}}

\def\epsi{\varepsilon}
\def\we{\wedge}

\def\de{\delta}

\def\part{\partial}

\def\R#1#2{ R^{#1}_{~~~#2} }
\def\PA#1#2{ P^{#1}_{A~~#2} }
\def\PS#1#2{ P^{#1}_{S~~#2} }

\def\Rp#1#2{ (R^+)^{#1}_{~~~#2} }

\def\Rm#1#2{ (R^-)^{#1}_{~~~#2} }
\def\Rinv#1#2{ (R^{-1})^{#1}_{~~~#2} }

\def\Rpm#1#2{(R^{\pm})^{#1}_{~~~#2} }

\def\Rb{{\bf \mbox{\boldmath $R$}}}
\def\Rbo{{\bf \mbox{\boldmath $R$}}}

\def\Rh{{\hat R}}

\def\Rhat#1#2{ \Rh^{#1}_{~~~#2} }

\def\L#1#2{ \La^{#1}_{~~~#2} }

\def\Rhatinv#1#2{ (\Rh^{-1})^{#1}_{~~~#2} }

\def\La{\Lambda}
\def\Rha{{\hat R}}
\def\ff#1#2#3{f_{#1~~~#3}^{~#2}}
\def\MM#1#2#3{M^{#1~~~#3}_{~#2}}

\def\cchi#1#2{\chi^{#1}_{~#2}}

\def\ome#1#2{\om_{#1}^{~#2}}

\def\RRhat#1#2#3#4#5#6#7#8{\La^{~#2~#4}_{#1~#3}|^{#5~#7}_{~#6~#8}}

\def\LL#1#2#3#4#5#6#7#8{\La^{~#2~#4}_{#1~#3}|^{#5~#7}_{~#6~#8}}

\def\Cb{\bf \mbox{\boldmath $C$}}
\def\CC#1#2#3#4#5#6{{\Cb}_{~#2~#4}^{#1~#3}|_{#5}^{~#6}}
\def\cc#1#2#3#4#5#6{C_{~#2~#4}^{#1~#3}|_{#5}^{~#6}}

\def\ZZ#1#2#3#4#5#6#7#8{Z^{~#2~#4}_{#1~#3}|^{#5~#7}_{~#6~#8}}
\def\C#1#2{ {\bf \mbox{\boldmath $C$}}_{#1}^{~~~#2} }
\def\c#1#2{ C_{#1}^{~~~#2} }

\def\DR{\Delta_R}
\def\DL{\Delta_L}
\def\f#1#2{ f^{#1}_{~~#2} }

\def\T#1#2{ T^{#1}_{~~#2} }

\def\M#1#2{ M_{#1}^{~#2} }

\def\rminus{r^{-1}}

\def\D{\Delta}

\def\Dp{\Delta^{\prime}}

\def\ep{\epsi^{\prime}}
\def\kp{{\kappa^{\prime}}}
\def\kpm{\kappa^{\prime -1}}

\def\km{\kappa^{-1}}

\def\rone{r \rightarrow 1}

\def\SOqrNt{SO_{q,r}(N+2)}

\def\ISOqrN{ISO_{q,r}(N)}
\def\ISpqrN{ISp_{q,r}(N)}

\def\SOqrN{SO_{q,r}(N)}

\def\SqrNt{S_{q,r}(N+2)}

\def\Dtwo{\Delta_{N+2}}
\def\epsitwo{\epsi_{N+2}}
\def\kappatwo{\kappa_{N+2}}

\def\Lpm#1#2{L^{\pm #1}_{~~~#2}}

\def\LLpm{L^{\pm}}

\def\Lp#1#2{L^{+ #1}_{~~~#2}}
\def\Lm#1#2{L^{- #1}_{~~~#2}}
\def\Dcal#1#2{{\cal D}^{#1}_{~#2}}

\def\chit{{\partial}}

\def\n2{{{N+1} \over 2}}
\def\ap{a^{\prime}}

\def\bu{\bullet}
\def\ci{\circ}
\def\square{{\,\lower0.9pt\vbox{\hrule \hbox{\vrule height 0.2 cm
\hskip 0.2 cm \vrule height 0.2 cm}\hrule}\,}}
\def\sma#1{\mbox{\footnotesize #1}}
\def\Q.E.D.{\rightline{$\Box$}}

\def\Nmezzi{{N\over 2}}
\def\fmi#1#2{f^{#1}_{~#2}}
\def\LM#1#2{\Lambda^{#1}_{~~~#2}}
\def\Acal{{\cal A}}
\def\parleft{{\stackrel{\leftarrow}{\chit}}}
\def\lati{{\tilde \la}}
\def\muti{{\tilde \mu}}
\def\le{\langle}
\def\re{\rangle}
\def\ISO{ISO_{q,r}(N)}

\begin{document}

\begin{titlepage}
\vskip -1cm
\rightline{DFTT-29/97}
\rightline{LBNL-40824}
%\rightline{April 1998}
\vskip 1em
\begin{center}
{\large\bf Quantum Orthogonal Planes:\\[-0.1em]
 {\boldmath $ISO_{q,r}(N)$} and {\boldmath $SO_{q,r}(N)$} --  
Bicovariant
Calculi\\[.2em]
and Differential Geometry on Quantum Minkowski Space}
\\[2em]
Paolo Aschieri \\[.4em]
{\sl Theoretical Physics Group, Physics Division\\
Lawrence Berkeley National Laboratory, 1 Cyclotron Road \\
Berkeley, California 94720, USA.}\\[1em]
Leonardo Castellani\\[.4em]
{\sl Dipartimento di Scienze e Tecnologie Avanzate*,
 Universit\`a di Torino; \\
Dipartimento di Fisica Teorica and Istituto Nazionale di
Fisica Nucleare\\
Via P. Giuria 1, 10125 Torino, Italy.} \\[1em]
Antonio Maria Scarfone\\[.4em]
{\sl Dipartimento di Fisica, Politecnico di Torino\\
Corso Duca degli Abruzzi 24,  10129 Torino, Italy}\\[2em]
\end{center}
\begin{abstract}
We construct differential calculi on multiparametric quantum
orthogonal planes in any dimension $N$. These calculi are
 bicovariant under the action of the full inhomogeneous
(multiparametric) quantum group $ISO_{q,r}(N)$, and do contain
dilatations.
 If we require bicovariance  only
under the quantum orthogonal group $SO_{q,r}(N)$, the calculus
on the $q$-plane can be expressed in terms of its coordinates $x^a$,
differentials $dx^a$ and partial derivatives $\partial_a$ without
the need of dilatations,  thus generalizing
known results to the multiparametric case.

Using real forms that lead to the signature
$(n+1,m)$ with $m\!=\!n-1,$ $n,n+1$, we find
$ISO_{q,r}(n+1, m)$ and $SO_{q,r}(n+1,m)$  bicovariant calculi
on the multiparametric quantum spaces. The particular case of
the quantum Minkowski space $ISO_{q,r}(3,1)/SO_{q,r}(3,1)$ is
treated
in detail.

The conjugated partial derivatives $\partial_a^*$ can
be expressed as linear combinations of the $\partial_a$.
This allows a deformation of the phase-space where no additional
operators
(besides $x^a$ and $p_a$) are needed.
\end{abstract}
\noi{April 1998}~~\hskip 10cm~q-alg/9709032

\centerline{{\small to appear on {\bf European Phys. Jou. C},  
Particles and
Fields }}
\vskip .2cm
\noi \hrule

\vskip .2cm
\noi{\small e-mail: aschieri@lbl.gov, castellani@to.infn.it,
scarfone@polito.it}

\noi$^*${\small \sl II Facolt\`a di Scienze M.F.N., sede di  
Alessandria}
\end{titlepage}
\newpage
\setcounter{page}{1}

%%%%%%%%%%%%%%%%%%%%%%%%%%%%%%
\sect{Introduction}
%%%%%%%%%%%%%%%%%%%%%%%%%%%%%%

Non commutativity of spacetime at the microscopic level could
provide
an
effective
regularization of gravity, in alternative to discretization methods.
It is
suggestive that a
non commutative structure of spacetime emerges in non-perturbative
attempts to
describe string theories \cite{Mtheorync}.

In this paper we use the non-commuting geometry \cite{noncom} of
quantum
groups
\cite{qgroups,FRT}, as defined by their differential calculi
\cite{Wor}
-\cite{Schupp}, to
derive the noncommuting differential geometry of the
multiparametric quantum
orthogonal planes in any dimension.
We then study real forms that are consistent  with the differential
calculus
and finally specialize our treatment to the multiparametric
quantum Minkowski space.

The necessary prerequisite for the work presented here  has been the
construction of inhomogeneous quantum groups of the orthogonal type
 $ISO_{q,r}(N)$ and of their corresponding bicovariant calculi.
This has been achieved in past publications
\cite{inhom1,inhom2,inhom3,inhom4} via a projection from the known
multiparametric
orthogonal groups $SO_{q,r}(N+2)$, and has provided an $R$ matrix
formulation
for the inhomogeneous case.

Other ref.s on inhomogeneous $q$-groups can be found in
\cite{inhom,tesi}.
For  multiparametric quantum groups see. ref.s
\cite{multiparam1}-\cite{multiparam2}.

In general, i.e. without any restrictions on the deformation
parameters,  inhomogeneous
groups of the orthogonal type contain dilatations. It is however
possible to  avoid
 dilatations if one fixes some of the parameters (including the $r$
parameter appearing in the off-diagonal terms of the $R$-matrix)
equal to one,
their classical
 value. The case $r=1$ corresponds to a
``quasi-classical" structure, for which the original braiding matrix
$\Rha$
becomes diagonal
(the corresponding deformations are then called {\sl twistings}).
In
this case
it is possible to construct a bicovariant calculus on $ISO_q(N)$,
and
consequently on
$q$-Minkowski space \cite{inhom2,inhom4}.

We present here a bicovariant calculus on the full multiparametric
$ISO_{q,r}(N)$
 {\sl without the restriction}  $r=1$. This calculus, however,  is
trivial on
the $SO_{q,r}(N)$
quantum subgroup: it can really be seen as a non-trivial calculus
only on the
coset $Fun_{q,r}[ISO(N)/SO(N)]$, i.e. on the quantum orthogonal
plane. For $r\not= 1$ this
$ISO_{q,r}(N)$--bicovariant calculus on the quantum plane
necessarily contains
dilatations.

If we require only $SO_{q,r}(N)$ bicovariance  [more precisely
right covariance under $ISO_{q,r}(N)$ and
left covariance only under $SO_{q,r}(N)$],
the calculus can be expressed in terms of coordinates $x$,
differentials $dx$ and partial derivatives $\partial$,
{\sl without the need of dilatations}.  In this case
the
$q$-commutations between $x$, $dx$ and $\partial$
close by themselves, and in
fact generalize to the multiparametric case  the known results
of ref.s
\cite{Manin,qplane,Fiore}.  Here these results emerge from the
broader setting of the bicovariant calculus on $ISO_{q,r}(N)$.
In this context we are able to explicitly relate the partial
derivatives $\partial$ to the $ISO_{q,r}(N)$ $q$-Lie algebra
generators.

It is natural to expect that a $*$-structure compatible with
$ISO_{q,r}(N)$
and with the bicovariant differential calculus will induce a well
behaved
$*$-conjugation on the differential calculus on the quantum plane,
acting linearly on the partial derivatives
$\partial_a, \;\sma{($a=1,...N$)}$.
The  conjugations that
give the real forms  $ISO_{q,r}(n+1,n-1)$,
$ISO_{q,r}(n,n)$ and  $ISO_{q,r}(n,n+1)$
are  consistent with the $ISO_q(N)$ bicovariant differential
calculus.
Using these conjugations one can define real coordinates $X^a$ and
hermitian partial derivative operators $P_a\sim \partial_a$, i.e.
momenta.
The conjugated $\partial^*_a$ are derived from the conjugation of
the
$ISO_{q,r}(N)$ $q$-Lie algebra generators and can be simply expressed
as
linear combinations of the $\partial_a$, {\sl without the need of
introducing an extra operator} as done in ref.s \cite{Lorek}.
The q-commutations of the momenta $P$ with the coordinates $X$
define a deformed phase-space that could be studied in the same
spirit as in ref.s \cite{Lorek}.

We will be concerned with the conjugation that gives the
$ISO_{q,r}(n-1,n+1)$ calculus and in particular induces a
differential calculus on the
$q$-Minkowski space. To retrieve the other conjugations, both for
$N$=even
and $N$=odd,
just take $\Dcal{A}{B}=\delta^A_B$ in the formulae where
$\Dcal{A}{B}$
appears.

In Section 2 we recall briefly the structure of the $ISO_{q,r}(N)$
quantum groups. Their
differential calculi are discussed in Section 3, and finally in
Sections 4 and 5 we present
the bicovariant calculi on quantum orthogonal planes in full detail.
In Appendix  A we
specialize our results to the four-dimensional quantum Minkowski
space, and list all the relevant formulas for its non-commuting
differential geometry.

%%%%%%%%%%%%%%%%%%%%%%%%%%%%%%%%%%%%%%%%%%%%%%%%
\sect{The quantum inhomogeneous group
$\ISOqrN$}
%%%%%%%%%%%%%%%%%%%%%%%%%%%%%%%%%%%%%%%%%%%%%%%%

{}~~~An $R$-matrix formulation for the
quantum inhomogeneous groups $\ISOqrN$ and $\ISpqrN$
was obtained in ref. \cite{inhom3}, in terms of
the $\R{AB}{CD}$  matrix for the
$SO_{q,r}(N+2)$ and $Sp_{q,r}(N+2)$ quantum groups.
We recall here the results for $SO_{q,r}(N+2)$. The quantum
inhomogeneous group
$\ISOqrN$  is freely generated by the
non-commuting matrix elements $\T{A}{B}$
 [{\small A}=$(\circ, a, \bullet)$, with $a=1,...N$)] and the
identity
$I$,
modulo the relations:
\eq
\T{a}{\circ}=\T{\bullet}{b}=\T{\bullet}{\circ}=0 ,
\label{Tprojected}
\en
\noi the $RTT$ relations
\eq
\R{AB}{EF} \T{E}{C} \T{F}{D} = \T{B}{F} \T{A}{E} \R{EF}{CD},
\label{RTTbig}
\en
\noi and the orthogonality relations
\eq
C^{BC} \T{A}{B}  \T{D}{C}= C^{AD} ,~~~
C_{AC} \T{A}{B}  \T{C}{D}=C_{BD} \label{CTTbig}
\en
The matrix $R$ controls the non-commutativity of the $\T{A}{B}$
elements,
and its entries depend continuously on a set of parameters $r,
q_{AB}$
($q_{AB}$ appearing only in the diagonal part of the $R$ matrix).
For
$r \rightarrow 1$, $q_{AB} \rightarrow 1$ (the ``classical limit"),
$\R{AB}{CD} \rightarrow \de^A_C \de^B_D$. The quantum metric
$C_{AB}$
and  its inverse $C^{AB}$ depend only on $r$ and  are
given in \cite{FRT}.

The co-structures of $\ISOqrN$ and $\ISpqrN$
are simply given by:
\eq
\D (\T{A}{B})=\T{A}{C} \otimes \T{C}{B},
{}~~\kappa (\T{A}{B})=C^{AC} \T{D}{C} C_{DB} ,~~
\epsi (\T{A}{B})=\de^A_B~. \label{costructuresbig}
\en
After decomposing the indices {\small A}=$(\circ, a, \bullet)$, and
defining:
\eq
u\equiv \T{\circ}{\circ},~~v\equiv
\T{\bullet}{\bullet},~~z\equiv
\T{\circ}{\bullet},~~
x^a \equiv \T{a}{\bullet},~~y_a \equiv \T{\circ}{a} \label{names}
\en
\noi the relations (\ref{RTTbig}) and (\ref{CTTbig}) become
\eqa
& &\R{ab}{ef} \T{e}{c} \T{f}{d} = \T{b}{f} \T{a}{e} \R{ef}{cd}
\label{PRTT11}\\
& &\T{a}{b} C^{bc} \T{d}{c}=C^{ad} I \label{PRTT31}\\
& &\T{a}{b} C_{ac} \T{c}{d} = C_{bd} I \label{PRTT32}
\ena
\eqa
& &\T{b}{d} x^a={r \over q_{d\bullet}} \R{ab}{ef} x^e \T{f}{d}
\label{PRTT33}\\
& &\PA{ab}{cd} x^c x^d=0 \label{PRTT13}\\
& &\T{b}{d} v={q_{b\bullet}\over q_{d\bullet}} v \T{b}{d}\\
& &x^b v=q_{b \bullet} v x^b \label{PRTT15}\\
& & uv=vu=I \label{PRTT21}\\
& &u x^b=q_{b\bullet} x^bu \label{PRTT22}\\
& &u \T{b}{d}={q_{b\bullet}\over q_{d\bullet}} \T{b}{d} u
\label{PRTT24}
\ena
\eq
y_b=-r^{\Nmezzi} \T{a}{b} C_{ac} x^c u \label{ipsilon}
\en
\eq
z=-{1\over {(r^{-\Nmezzi}+ r^{\Nmezzi-2})}} x^b C_{ba} x^a u
\label{PRTT44}
\en
\noi where $q_{a\bullet}$ are $N$  complex parameters
related by $q_{a\bullet} = r^2 /q_{\ap\bullet}$, with $\ap = N+1-a$.
The matrix $P_A$ in eq. (\ref{PRTT13}) is the $q$-antisymmetrizer
for
the
$B,C,D$ $q$-groups given by (cf.  (\ref{projBCD})):
\eq
\PA{ab}{cd}=- {1 \over {r+\rminus}}
(\Rhat{ab}{cd}-r\de^a_c \de^b_d + {r-r^{-1} \over
 r^{N-2} +1} C^{ab} C_{cd}). \label{PA}
\en
The last two relations (\ref{ipsilon}),  (\ref{PRTT44})
are constraints, showing that
the $\T{A}{B}$ matrix elements in eq. (\ref{RTTbig})
are really a {\sl redundant} set. This redundance
is necessary if we want to express the $q$-commuations
of the $\ISOqrN$ basic group elements
as $RTT=TTR$  (i.e. if we want an $R$-matrix
formulation).  We can take as independent generators the
elements
\eq
{}~~\T{a}{b} , x^a , v , u\equiv v^{-1}
\mbox{ and the identity }  I~~~~~(a=1,...N)
\en
The co-structures on the $\ISOqrN$ generators can be read from
(\ref{costructuresbig})
after decomposing the indices \sma{$A=\circ,a,\bu$}:
\eqa
& &\D (\T{a}{b})=\T{a}{c} \otimes \T{c}{b}~,~~
\D (x^a)=\T{a}{c} \otimes x^c + x^a \otimes v ~,
\label{Pcoproduct2}\\
& &\D (v)=v \otimes v~,~~\D (u)=u \otimes u~,
\ena
\eqa
& &\kappa (\T{a}{b})
=C^{ac} \T{d}{c} C_{db}~,\\
& &\kappa (x^a)
=- \kappa (\T{a}{c}) x^c u~,~~\kappa (v)= u~,~~\kappa (u)= v~,
\ena
\eq
\epsi (\T{a}{b})=\de^a_b ~,~~\epsi (x^a)=0 ~,~~
\epsi (u)=\epsi (v)=\epsi(I)=1 ~.\label{cfin}
\en
\sk
In the commutative limit $q\rightarrow 1 , r\rightarrow 1$ we
recover the algebra of functions
on $ISO(N)$ (plus the dilatation $v$ that can be set to
the identity).
\sk
\noi {\bf Note 2.1 : }  as shown in ref. \cite{inhom3},
the quantum group $\ISOqrN$
can be derived as the quotient
\eq
\frac{\SOqrNt}{H}~,
\label{quotient}
\en
\noi where $H$ is the Hopf ideal in
 $\SOqrNt$ of
all sums of monomials containing at least an element of the kind
$\T{a}{\circ}, \T{\bullet}{b}, \T{\bullet}{\circ}$.
The Hopf structure of the groups in the numerators of
(\ref{quotient})
is naturally inherited by the quotient groups \cite{Sweedler}.

We denote by $P$ the canonical projection
\eq
P ~:~~  \SqrNt\longrightarrow \SqrNt/{H}
\en
\noi It is a Hopf algebra epimorphism because $H=Ker(P)$ is a Hopf
ideal. Then
any element of ${\SqrNt/H}$ is of the form $P(a)$ and
\eq
P(a)+P(b)\equiv P(a+b) ~;~~ P(a)P(b)\equiv P(ab) ~;~~
\mu P(a)\equiv P(\mu a),~~~\mu \in \mbox{\bf C} \label{isoalgebra}
\en
\eq
\D (P(a))\equiv (P\otimes P)\Dtwo(a) ~;~~ \epsi(P(a))
\equiv\epsitwo(a) ~;~~
\kappa(P(a))\equiv P(\kappatwo(a)) \label{co-iso}
\en
\noi where we  indicate by
$\Dtwo$, $\epsitwo$ and $\kappatwo$ the
co-structures of  $\SOqrNt$.
Eq.s (\ref{PRTT11}) - (\ref{PRTT44})  have been obtained in
\cite{inhom3}
by taking the $P$
projection of the $RTT$ and $CTT$ relations of
$\SqrNt$, with the notation $u\equiv P(\T{\circ}{\circ}),~v\equiv
P(\T{\bullet}{\bullet}),~z\equiv P(\T{\circ}
{\bullet}), ~x^a \equiv P(\T{a}{\bullet}),~y_a \equiv
P(\T{\circ}{a}),
{}~T^a{}_b\equiv P(T^a{}_b) ~;~~ I\equiv P(I) ~;~ 0\equiv P(0)$,
cf. (\ref{names}).
\sk
\noi {\bf Note 2.2: } From the commutations
(\ref{PRTT22}) - (\ref{PRTT24})
 we see that
one can set $u=I$ only when $q_{a\bullet}=1$ for all $a$.
{}From $q_{a\bullet} = r^2 /q_{\ap\bullet}$
this implies also $r=1$.
\sk
\noi {\bf Note 2.3:} Eq.s (\ref{PRTT13}) are just the
multiparametric
orthogonal quantum plane commutations. Choosing as free indices
$\left({}^{\bu b}_{~\bu\bu}\right)$ in  (\ref{RTTbig}) yields
 $zx^b=q_{b\bu}x^bz$ and therefore  the (squared) length
element $L=x^aC_{ab}x^b$ commutes with the $x$ elements.
Similarly  we find  $L\T{a}{d}=(q_{d\bu}/r)^2\T{a}{d}L$ and
$Lu=r^{-2}uL,~Lv=r^2vL$.
\sk
\noi {\bf Note 2.4:} Two conjugations
(i.e. algebra antihomomorphisms, coalgebra
homomorphisms and involutions, satisfying $\kappa(\kappa(T^*)^*)=T$)
exist on $ISO_{q,r} (N)$ , inherited from the corresponding ones on
$SO_{q,r}(N+2)$
\cite{inhom3,inhom4}. We recall here their action on the group
generators
$\T{A}{B}$ and the corresponding restrictions on the parameters:
\sk
$\bullet$~~trivially as $T^*=T$; corresponds  to the real forms
$ISO_{q,r}(n,n;\Rb)$ and $ISO_{q,r}(n,n+1;\Rb)$. Compatibility  with
the
$RTT$ relations (\ref{RTTbig}) requires
$|q_{ab}|=|q_{a\bullet}|=|r|=1$.
\sk
$\bullet$~~Only for $N=2n$ even: $ (\T{A}{B})^*={\cal D}^A_{~C}
\T{C}{D} {\cal D}^D_{~B} $,  ${\cal D}$ being the matrix that
exchanges the index $n$ with the index $n+1$; extends to the
inhomogeneous multiparametric case the conjugation
proposed by the authors of ref. \cite{Firenze1} for $SO_{q}(2n,{\bf
C})$, and corresponds to the real form
$ISO_{q,r}(n+1,n-1;\Rb)$.

 Explicitly: $(\T{a}{b})^*={\cal D}^a_{~c} \T{c}{d} {\cal D}^d_{~b}
,
 (x^a)^*=\Dcal{a}{b} x^b, u^*=u, v^*=v, z^*=z$.
Compatibility with the $RTT$  relations (\ref{RTTbig})
requires:
\eq
({\bar R})_{n \leftrightarrow n+1}=R^{-1},~~~\mbox{ i.e. }~~{\cal
D}_1 {\cal
D}_2 R_{12}
{\cal D}_1 {\cal D}_2 = {\overline R_{12}^{-1} } \label{Rprop2}
\en
which implies $|r|=1$;
$|q_{a\bullet}|=1$ for $a \not= n, n+1$;
$|q_{ab}|=1$ for $a$
and $b$
both different from $n$ or $n+1$;  $q_{ab}/r \in {\Rb}$
when at least one of the indices $a,b$ is equal
to $n$ or $n+1$;   $q_{a\bullet}/r \in \Rb$ for $a=n$ or $a=n+1$.
Compatibility with the $CTT$ relations (\ref{CTTbig}) is ensured
by ${\cal D}C{\cal D}$ and ${\bar C}=C^T$ (due to $|r|=1$).

In particular,  the quantum Poincar\'e group $ISO_{q,r}(3,1;\Rb)$
is obtained by setting $|q_{1\bullet}|=|r|=1$,
$q_{2\bullet}/r \in \Rb$, $q_{12}/r  \in \Rb$.

According to Note 2.2,  a {\sl dilatation-free} $q$-Poincar\'e
group
is found after the further restrictions
$q_{1\bullet}=q_{2\bullet}=r=1$.
The only free parameter remaining is then $q_{12} \in \Rb$.
\sk
%%%%%%%%%%%%%%%%%%%%%%%%%%%%%%%%%%%%%%%%%%%%%%
\sect{Bicovariant calculus on simple $q$-groups}
%%%%%%%%%%%%%%%%%%%%%%%%%%%%%%%%%%%%%%%%%%%%%%
The bicovariant differential calculus on the
$q$-groups of the $A,B,C,D$ series can be formulated in terms
of the corresponding $R$-matrix, or equivalently in terms of
the $\LLpm$ functionals defined by:
\eq
\Lpm{A}{B} (\T{C}{D})= \Rpm{AC}{BD}~, ~~~
\Lpm{A}{B} (I)=\de^A_B \label{LonT}
\en
\noi with \footnote{for the $B,C, D$ series. For the
$q$-groups of the $A$ series there is more freedom in choosing $R^+$
and $R^-$,
see for ex. \cite{inhom1}}
\eq
\Rp{AC}{BD} \equiv \R{CA}{DB} ~,~~~
\Rm{AC}{BD} \equiv \Rinv{AC}{BD}~. \label{Rplusminus}
\en
To extend the definition (\ref{LonT})
to the whole Hopf algebra $\Acal$ we set
\eq
\Lpm{A}{B} (ab)=\Lpm{A}{C} (a) \Lpm{C}{B} (b)
{}~~~\forall a,b\in  \Acal  ~.
\label{Lab}
\en
These functionals generate the Hopf algebra $\Acal'$ paired to
$\Acal$, with
$\D'(\Lpm{A}{B})=\Lpm{A}{C} \otimes \Lpm{C}{B}~,$
$\epsi'(\Lpm{A}{B})=\delta^A_B$ and $\kappa'(\LLpm)=(\LLpm)^{-1}$.
\sk
We briefly recall how to construct a bicovariant
calculus.  The general procedure can be found in ref.s
\cite{Jurco,Schupp}, or, in the notations we adopt
here, in ref. \cite{Aschieri1}.
It realizes the axiomatic construction of ref. \cite{Wor}.
\sk
The functionals
\eq
\ff{A_1}{A_2B_1}{B_2} \equiv \kappa (\Lp{B_1}{A_1}) \Lm{A_2}{B_2}.
\label{defff}
\en
and the elements of $\Acal$:
\eq
\MM{B_1}{B_2A_1}{A_2} \equiv \T{B_1}{A_1} \kappa (\T{A_2}{B_2}).
\label{defMM}
\en
satisfy the following relations, called {\sl bicovariant bimodule
conditions},
where for simplicity
we use the adjoint indices $ i,j,k,...$
with ${}^i={}_A^{~B},~{}_i={}^A_{~B}~:$
\eqa
& & \D' (\f{i}{j}) = \f{i}{k} \otimes \f{k}{j}  \label{propf1}~~;~~~
\epsi' (\f{i}{j}) = \del^i_j~,  \\
& & \Delta (\M{j}{i}) = \M{j}{l} \otimes \M{l}{i} \label{copM}~~;~~~
\epsi (\M{j}{i}) = \delta^i_j~, \\
& & \M{i}{j} (a * \f{i}{k})=(\f{j}{i} * a) \M{k}{i} \label{propM}
\ena
\noi The star product between a functional on $\Acal$ and an
element of $\Acal$ is defined as:
\eq
 \chi * a \equiv (id \otimes \chi) \D (a),~~~
a * \chi \equiv (\chi \otimes id) \D (a),~~~
a \in \Acal, ~\chi \in \Acal '
\en
Relation  (\ref{propM}) is easily checked for
$a=\T{A}{B}$ since in
this case it is implied by the $RTT$ relations;
it holds for a generic $a$ because of  property
(\ref{propf1}).
\sk
The space of {\bf  quantum one-forms} is defined as a left
$\Acal$-module $\Ga$
freely generated  by the  {\sl symbols} $\ome{A_1}{A_2}$:
\eq
\Ga\equiv \{a^{A_1}_{~A_2}\ome{A_1}{A_2}\}~ ,~~~a^{A_1}_{~A_2} \in
\Acal \label{symbols}
\en
{\bf Theorem 3.1} (due to Woronowicz, see last ref. in \cite{Wor}):
because of
the
properties  (\ref{propf1}), $\Ga$
becomes a bimodule over $\Acal$ with the following
right multiplication:
\eq
\ome{A_1}{A_2} a=(\ff{A_1}{A_2B_1}{B_2} * a)~\ome{B_1}{B_2},
\label{omea}
\en
in particular:
\eq
\ome{A_1}{A_2} \T{R}{S}= \Rinv{TB_1}{CA_1}
\Rinv{A_2C}{B_2S} \T{R}{T}\ome{B_1}{B_2} \label{commomT}
\en
Moreover, because of properties (\ref{copM}) and (\ref{propM}),
 we can  define a left and a right action of
$\Acal$ on $\Ga$:
\eqa
{}\!\!\!\!\!\!
& & \DL:\Ga\rightarrow \Acal\otimes
\Ga~~;~~~\DL(a \ome{A_1}{A_2} b)
\equiv \D(a)\,(I\otimes \ome{A_1}{A_2})\,\D(b)~,
\label{deltaLomega} \\
{}\!\!\!\!\!\!& & \DR:\Ga\rightarrow  \Ga\otimes \Acal
{}~~;~~~\DR(a \ome{A_1}{A_2} b)
\equiv \D(a)\,(\ome{B_1}{B_2}\otimes
\MM{B_1}{B_2A_1}{A_2})\,\D(b)~.
\label{deltaRomega}
\ena
These actions commute, i.e. $(id \otimes \DR) \DL = (\DL \otimes id)
\DR$,
and give a bicovariant bimodule structure to $\Ga$.
[Property (\ref{propM}) is a sufficient and necessary condition for
$\DR(\rho b)=\DR(\rho)\D(b)$].
\sk
The elements $\M{j}{i}$ can be used to
build a right-invariant basis of $\Ga$. Indeed the $\eta^i$ defined
by
\eq
\eta^i \equiv \kappa^{-1}(\M{j}{i})\om^j  \label{eta}
\en
are right invariant: $\DR(\eta^i)=\eta^i\otimes I$
[use $\kappa^{-1}(a_2)a_1=\epsi(a)$];
moreover every element of $\Ga$ can be  written
as  $\rho = \eta^ib_i$
or $a_i\eta^i$ where
$b_i$ and $a_i$ are uniquely determined.
It can be shown that the functionals $\f{i}{j}$
satisfy:
\eqa
& & \eta^i b = (b * \f{i}{j}) \eta^j\label{etab}\\
& & a \eta^i = \eta^j [a*(\f{i}{j} \circ \kappa)],\label{etabbb}
\ena
\noi where $a*f \equiv (f \otimes id) \D(a)$.
\sk
The {\bf exterior derivative}
$d ~:~ \Acal \longrightarrow \Ga$ can be defined via the element
$\tau\equiv \sum_A \ome{A}{A} \in \Ga$. This element
is
easily shown to be
left and right-invariant:
\eq
\DL(\tau)=I\otimes\tau~~;~~~
\DR(\tau)=\tau\otimes I
\en
and defines the derivative  $d$  by
\eq
da=\rinv [\tau a - a \tau]. \label{defd1}
\en
\noi The factor $\rinv$ is necessary for a correct
classical limit $\rone$.
It is immediate to prove the Leibniz rule
\eq
d(ab)=(da)b+a(db),~~\forall a,b\in \Acal ~. \label{Leibniz}
\en
Two other expressions for the derivative are given by:
\eq
da = ( \cchi{A_1}{A_2} * a) ~\ome{A_1}{A_2} \label{defd2}~,
\en
\eq
da = -\eta_{A_1}^{~A_2}(a* \kappa'(\cchi{A_1}{A_2})\,)
\label{defd2eta}~
\en
where the linearly independent elements
\eq
\cchi{A}{B} = \rinv [\ff{C}{CA}{B}-\de^A_B \epsi]
\label{defchi}
\en
are the tangent vectors such that the left-invariant vector fields
$\chi^A{}_B *$
are dual to the left-invariant
one-forms $\ome{A_1}{A_2}$, and similarly the right-invariant vector
fields
$*\kappa'(\cchi{A_1}{A_2})$ are dual to the right-invariant one forms
$-\eta_{A_1}^{~A_2}$.
The equivalence of  (\ref{defd1}) and (\ref{defd2}) can be
shown by using the rule (\ref{omea}) for $\tau a$ in
the right-hand side of (\ref{defd1}). Expression (\ref{defd2eta})
is related to (\ref{defd2}) via
$\chi_i*a=(a*\chi_j)\kappa^{-1}(M_i{}^j)$,
eq.s (\ref{eta}),
(\ref{etabbb}) and $\kp (\chi_i)=-\chi_j \kp (\f{j}{i})$.

\noi Using (\ref{defd2}) we can compute the exterior derivative
on the basis elements $\T{A}{B}$:
\eq
d ~\T{A}{B}=\rinv [\Rinv{CR}{ET} \Rinv{TE}{SB}
\T{A}{C}-\de^R_S \T{A}{B}]
{}~\ome{R}{S}  \label{dTAB}
\en
\sk
Every element  $\rho$ of $\Gamma$, which by definition is
written in a unique way as
$\rho=a^{A_1}{}_{A_2}\ome{A_1}{A_2}$,  can also be written as
\eq
\rho=\sum_k a_k db_k  \label{propd}
\en
\noi for some $a_k,b_k$ belonging to $\Acal$. This can
be proven directly by inverting
the relation (\ref{dTAB}).
\sk
Due to the bi-invariance of $\tau$ the derivative operator $d$ is
compatible with the actions $\DL$ and $\DR$:
\eq
\DL (adb)=\D(a)(id\otimes d)\D(b)~~,~~~\DR (adb)=\D(a)(d\otimes
id)\D(b)~,\label{dbicov}
\en
these two properties express the fact that $d$
commutes with the left
and right action of the quantum group, as in the classical case.
\sk
\noi {\bf Remark :} The properties (\ref{Leibniz}),
(\ref{propd})  and (\ref{dbicov}) of
the exterior derivative (\ref{defd2})
realize the axioms of a first-order bicovariant differential
calculus \cite{Wor}.
\sk
The {\bf tensor product} between elements $\rho,\rhop \in \Ga$
is defined to
have the properties $\rho a\otimes \rhop=\rho
\otimes a \rhop$, $a(\rho
\otimes \rhop)=(a\rho) \otimes \rhop$ and
$(\rho \otimes \rhop)a=\rho
\otimes (\rhop a)$. Left and right actions on
$\Ga \otimes \Ga$ are
defined by:
\eqa
{}\!\!\!\!\!\!\!
&&\DL (\rho \otimes \rhop)\equiv \rho_1 \rhop_1
\otimes \rho_2 \otimes
\rhop_2,~~~\DL: \Ga \otimes \Ga \rightarrow \Acal
\otimes\Ga\otimes\Ga
\label{DLGaGa}\\
{}\!\!\!\!\!\!\!
&&\DR (\rho \otimes \rhop)\equiv \rho_1 \otimes \rhop_1
 \otimes \rho_2
\rhop_2,~~~\DR: \Ga \otimes \Ga \rightarrow
\Ga\otimes\Ga\otimes \Acal
\label{DRGaGa}
\ena
\noi where  $\rho_1$, $\rho_2$, etc., are defined by:
$$
\DL (\rho) = \rho_1 \otimes \rho_2,~~~\rho_1\in \Acal,~\rho_2\in
\Ga~~;~~~
\DR (\rho) = \rho_1 \otimes \rho_2, ~\rho_1\in \Ga,~\rho_2\in
\Acal~.$$
The extension to $\Ga^{\otimes n}$ is straightforward.
\sk
The  {\bf exterior product} of one-forms is consistently defined as:
\eq
\ome{A_1}{A_2} \we \ome{D_1}{D_2}
\equiv \ome{A_1}{A_2} \otimes \ome{D_1}{D_2}
- \RRhat{A_1}{A_2}{D_1}{D_2}{C_1}{C_2}{B_1}{B_2}
\ome{C_1}{C_2} \otimes \ome{B_1}{B_2} \label{exteriorproduct}
\en
\noi where the $\Lambda$ tensor is given by:
\eq
\begin{array}{rl}
\LL{A_1}{A_2}{D_1}{D_2}{C_1}{C_2}{B_1}{B_2}
\equiv & \ff{A_1}{A_2B_1}{B_2} (\MM{C_1}{C_2D_1}{D_2}) = \\
=&
d^{F_2} d^{-1}_{C_2} \R{F_2B_1}{C_2G_1} \Rinv{C_1G_1}{E_1A_1}
    \Rinv{A_2E_1}{G_2D_1} \R{G_2D_2}{B_2F_2} \label{Lambda}
\end{array}
\en
\noi $d^A$ being the entries of the diagonal matrix
$D^A_{~B} \equiv C^{AC} C_{BC}$.  From the formula
(\ref{exteriorproduct})
we can find the $q$-commutations (generalizing the classical
anticommutations)
of products of one-forms $\om$ in terms of a ``flip" operator (see
the second
ref. in \cite{qgtorino}):
\eq
\om^i \we \om^j= - Z^{ij}_{~~kl} \om^k \we \om^{l}  \label{flip}
\en
\sk
Using the exterior product we can define the {\bf exterior
differential on $\Ga$} :
\eq
d~:~\Gamma \rightarrow \Gamma \we \Gamma
{}~~~;~~~~d (a_k db_k) = da_k \we db_k \label{defdonga}
\en
\noi which can easily be extended to
$\Gamma^{\we n}$ ($d: \Gamma^{\we n}
\rightarrow \Gamma^{\we (n+1)}$, $\Gamma^{\we n}$ being
defined as in the
classical case but with the quantum permutation
(braid) operator $\La$ \cite{Wor}). The definition (\ref{defdonga})
is equivalent to the following :
\eq
d\theta = \rinv [\tau \we \theta - (-1)^k \theta \we \tau],
\label{defdgen}
\en
\noi where $\theta \in \Ga^{\we k}$. The exterior differential
 has the following properties:
\eq
d(\theta \we \thetap)=d\theta \we \thetap +
(-1)^k \theta \we d\thetap~~~;~~~~
d(d\theta)=0~,\label{propd2}
\en
\eq
\DL (\theta d\theta')=\DL(\theta)(id\otimes
d)\DL(\theta')\label{propd3}~~~;~~~~
\DR (\theta d\theta')=\DR(\theta)(d\otimes
id)\DR(\theta'),\label{propd4}
\en
\noi where $\theta \in \Ga^{\we k}$, $\thetap \in \Ga^{\we n}$.
\sk
The {\bf \mbox{\boldmath $q$}
-Cartan-Maurer equations} are found by
using
(\ref{defdgen}) in computing $d\ome{C_1}{C_2}$:
\eq
d\ome{C_1}{C_2}= \rinv (\ome{B}{B} \we \ome{C_1}{C_2} +
 \ome{C_1}{C_2} \we
\ome{B}{B}) \equiv
-\onehalf \cc{A_1}{A_2}{B_1}{B_2}{C_1}{C_2}
{}~\ome{A_1}{A_2} \we \ome{B_1}{B_2} \label{CartanMaurer}
\en
\noi with:
\eq
\cc{A_1}{A_2}{B_1}{B_2}{C_1}{C_2}= -{2\over (r-r^{-1})}
[ \ZZ{B}{B}{C_1}{C_2}{A_1}{A_2}{B_1}{B_2} + \de^{A_1}_{C_1}
\de^{C_2}_{A_2} \de^{B_1}_{B_2} ]\label{cc}
\en
To derive this formula we have used the flip operator $Z$
on $\ome{B}{B} \we \ome{C_1}{C_2}$.
\sk
Finally, we recall that the $\chi$ operators close
on the {\bf $q$-Lie algebra} :
\eq
\chi_i \chi_j - \L{kl}{ij} \chi_k \chi_l = \C{ij}{k} \chi_k
\label{bico1}
\en
\noi where the $q$-structure constants are given by
\eq
{\C{jk}{i}=\chi_k(\M{j}{i})} ~~\mbox{ explicitly : }~~
{\CC{A_1}{A_2}{B_1}{B_2}{C_1}{C_2}} =\rinv [- \de^{B_1}_{B_2}
\de^{A_1}_{C_1} \de^{C_2}_{A_2} +
\LL{B}{B}{C_1}{C_2}{A_1}{A_2}{B_1}{B_2}]. \label{CC}
\en
The $C$ structure constants
appearing in the Cartan-Maurer equations are in general related
to the $\Cb$ constants of the $q$-Lie algebra \cite{Aschieri1}:
\eq
\C{jk}{i}=\onehalf [\c{jk}{i}-\L{rs}{jk} \c{rs}{i}]~.
\en
\sk
Using the definitions (\ref{defchi}) and (\ref{defff}) it is
not difficult to find the co-structures
on the functionals $\chi$ and $f$:
\eq
\begin{array}{lcl}
\Dp (\chi_i)=\chi_j
\otimes \f{j}{i} + \epsi \otimes \chi_i
&~;~~ &\Dp (\f{i}{j})=\f{i}{k} \otimes \f{k}{j}~, \\
\ep(\chi_i)=0
&~;~~ & \ep (\f{i}{j}) = \del^i_j~, \\
\kp(\chi_i)=-\chi_j \kp(\f{j}{i})
&~;~~ & \kp (\f{k}{j}) \f{j}{i}= \de^k_i \epsi =
\f{k}{j} \kp (\f{j}{i})~.
\end{array}\label{copf}
\en
\noi Note
that in the $r,q \rightarrow 1$ limit $\f{i}{j}
\rightarrow \de^i_j \epsi$, i.e. $\f{i}{j}$ becomes
proportional to the
identity functional and formula (\ref{omea}),
becomes trivial, e.g. $\om^i a = a\om^i$ [use $\epsi * a=a$].
\sk
The  $*$-conjugation on $\Acal$ is canonically extended to a
conjugation
on the  Hopf algebra $\Acal'$ generated by the $\LLpm$ functional
and
paired to $\Acal$. The relation is
\eq
\phi^*(a)\equiv \overline{\phi(\kappa^{-1}(a^*))}\label{stardef}
\en
where the overline denotes complex conjugation.
Explicitly, on the $\chi$ functionals it reads \footnote{We
correct here a misprint of \cite{inhom4} where the factor
$r^{-N+1}$ instead of $r^{-N-1}$ appears.
Notice that there we have used the opposite convention
$\phi^*(a)\equiv \overline{\phi(\kappa(a^*))}$
instead of (\ref{stardef}).}  \cite{inhom4}
\eq
(\cchi{A}{B})^*=- r^{-N-1} \cchi{C}{D} \Dcal{F}{~B}
 \Dcal{A}{~G} \R{EG}{FC} D^D_{~E}~~~
\sma{\mbox{ for } $SO_{q,r}(n+2,n ;\Rbo)~~;~~~  2n+2=N+2$}
\label{conjchiAB}
\en
\noi with $D^D_{~E} \equiv C^{DF} C_{EF}$.
Since the conjugation is a linear operation on the $q$-Lie algebra,
it can be
extended via (\ref{defd2}) to the differential calculus as well
\cite{Wor}, \cite{tesi}. The unique antilinear involution $*$ on
$\Ga$
satisfies:
\eqa
& & (a\rho)^*=\rho^* a^*, ~(\rho a)^*=a^* \rho^*, (da)^*=d(a^*); \\
& & \DL (\rho^*)=\DL (\rho)^*,~~\DR (\rho^*)=\DR (\rho)^*
\ena
\noi where $(a\otimes \rho)^*=a^* \otimes \rho^*$ and $(\rho \otimes
a)^*=\rho^* \otimes a^*$.
It easily extends to
$\Gamma^{\we n}$, for ex. $(d\theta)^*=d \theta^*$ etc.
Inverting formulae (\ref{dTAB}) one can also find
the induced conjugation on the
left-invariant one-forms \cite{inhom4}.
The explicit relation between the $*$-structures on the $\chi$ and on
the
$\omega$
can be given using the duality
$\le \ome{A}{B},\cchi{C}{D}\re=\delta_A^C\delta_D^B$
between left-invariant vector fields and left-invariant one-forms
\cite{tesi}:
\eq
\le {\ome{A}{B}}^{\,*},\cchi{C}{D}\re=
-\overline{\le {\ome{A}{B}},{\cchi{C}{D}}^*}\re~.
\label{om*chi}\en

%%%%%%%%%%%%%%%%%%%%%%%%%%%%%%%%%%%%%%
\sect{Bicovariant  differential calculus on $\ISOqrN$}
%%%%%%%%%%%%%%%%%%%%%%%%%%%%%%%%%%%%%%
The existence of this calculus is simply due to the existence
of a subset of the functionals (\ref{defff}), vanishing on the Hopf
ideal $H$
(see Section 2) ,  and $M$ elements that
satisfy the bicovariant bimodule conditions
(\ref{propf1})-(\ref{propM}).
These are:
\eq
\ff{\bu}{\ci \bu}{\ci}~,~~\ff{\bu}{a \bu}{\ci}~,~~\ff{\bu}{\bu
\bu}{\ci}~,~~
\ff{\bu}{a \bu}{b}~,~~\ff{\bu}{\bu \bu}{b}~,~~\ff{\bu}{\bu
\bu}{\bu}~
\label{fin}
\en
\eq
\begin{array}{lll}
P(\MM{\bu}{\ci\bu}{\ci})=v^2   &P(\MM{\bu}{\ci\bu}{d})=0
&P(\MM{\bu}{\ci\bu}{\bu})=0     \\
P(\MM{\bu}{b \bu}{\ci})=vr^{-{N\over 2}}x^eC_{eb}      &P(\MM{\bu}{b
\bu}{d})=v\kappa(\T{d}{b})    &P(\MM{\bu}{b \bu}{\bu})=0     \\
P(\MM{\bu}{\bu\bu}{\ci})=-{1\over{r^N(r^{N\over 2}+r^{-{N\over
2}+2})}}x^eC_{ef}x^f
&P(\MM{\bu}{\bu\bu}{d})=v\kappa(x^d)     &P(\MM{\bu}{\bu\bu}{\bu})=I
{}~~~
\label{Min}
\end{array}
\en
Notice that only the couples of indices $\sma{$(\bu \ci),~(\bu b)$}$
and
$\sma{$(\bu\bu)$}$ appear
in (\ref{fin})-(\ref{Min}): these are therefore the only indices
involved
in the projected differential calculus on $ISO_{q,r}(N)$.
\sk
\noi {\bf Theorem 4.1:} the functionals $\f{i}{j}$ in (\ref{fin})
 annihilate the Hopf ideal $H$.
\sk
\noi {\sl Proof : } one first checks directly that the
functionals (\ref{fin}) vanish
on the generators ${\cal T}$ of the ideal  $H$, i.e. on ${\cal T} =
\T{a}{\circ}, \T{\bu}{b}, \T{\bu}{\ci}$.
This extends to any element of
$H=Ker(P)$, because of the property (\ref{propf1}).  Q.E.D.
\sk
 These functionals are then
well defined on the quotient $\ISOqrN = \SOqrNt / Ker(H)$,
in the sense that the ``projected" functionals
\eq
\fb : \ISOqrN \rightarrow {\Cb},~~
\fb (P(a)) \equiv f(a)~ , ~~~\forall a \in \SOqrNt \label{deffb}
\en
are well defined.
Indeed if $P(a)=P(b)$, then $f(a)=f(b)$ because
$f(Ker(P))=0$. This holds for any
 functional $f$ vanishing on $Ker(P)$.
\sk
The product $fg$ of  two generic functionals
vanishing on $KerP$ also vanishes
on $KerP$, because $KerP$ is a co-ideal (see ref. \cite{inhom3}):
 $fg(KerP)=(f\otimes g)\D_{N+2} (KerP)=0$.
Therefore $\overline{fg}$ is well defined on $\ISOqrN$; moreover,
[use
(\ref{co-iso})]
$\overline{fg}[P(a)]\equiv fg(a)= (\bar{f}\otimes \bar{g}) \D(P(a))
\equiv
 \bar{f} \bar{g}[P(a)] $,
\noi and the product of $\bar{f}$ and $\bar{g}$ involves
the coproduct $\D$ of $\ISOqrN$.
\sk
There is a natural way to introduce a coproduct on the $\fb$'s :
\eq
\D\fb[P(a)\otimes P(b)]\equiv\fb[P(a)P(b)]=\fb[P(ab)]
=f(ab)=\D_{N+2}f(a\otimes b) ~.
\en
It is also easy to show that
\eq
\D(\fb^i{}_j) = \fb^i{}_k\otimes \fb^k{}_j ~~\mbox{ i.e. }~
\fb^i{}_j[P(a)P(b)] = \fb^i{}_k[P(a)]\fb^k{}_j[P(b)]
\label{cofb}
\en
with $i,j,k$ running over the restricted set of indices
{\small $\bu
b,\bu\bu,\bu\ci$ }.
Indeed due to the vanishing of some $f$'s (a consequence of
upper and lower triangularity of $L^+$ and
$L^-$ respectively), formulae (\ref{copf}) and (\ref{propf1})
involve only the $f^i{}_j$ listed in (\ref{fin}). Then
\eq
\fb^i{}_j[P(a)P(b)]=\fb^i{}_j[P(ab)]=f^i{}_j(ab)
=f^i{}_k(a)f^k{}_j(b)=\fb^i{}_k[P(a)]\fb^k{}_j[P(b)]
\en
and (\ref{cofb}) is proved.
\sk
With abuse of notations we
will simply write $f$ instead of $\fb$. Then the
$f$  in (\ref{fin}) will be seen as functionals on
$\ISOqrN$.
\sk
\noi {\bf Theorem 4.2:} the left $\Acal$-module ($\Acal=\ISOqrN$)
 $\Ga$ freely generated
by $\om^i \equiv \ome{\bu}{a}, \ome{\bu}{\bu}$ and $\ome{\bu}{\ci}$
is a bicovariant bimodule over $\ISOqrN$ with right multiplication:
\eq
\om^i a = (\f{i}{j} * a) \om^j,~~a \in \ISOqrN \label{omia}
\en
\noi where the $\f{i}{j}$ are given in (\ref{fin}), the $*$-product
is computed with the co-product $\D$ of $\ISOqrN$,  and
the left and right actions
of $\ISOqrN$ on $\Ga$ are given by
\eqa
& &\DL (a_i \om^i) \equiv \D(a_i) I \otimes \om^i
 \label{DLin}\\
& &\DR (a_i \om^i) \equiv \D(a_i) \om^j \otimes \M{j}{i}
\label{DRin}
\ena
\noi where the $\M{j}{i}$ are given in (\ref{Min}).
\sk
{\sl Proof:}
by showing  that  the
functionals $f$
and the elements $M$ listed in (\ref{fin}) and
(\ref{Min})
satisfy the properties (\ref{propf1})-(\ref{propM})
(cf.  Theorem  3.1).
 It  is straightforward to verify directly that
 the elements $M$
in (\ref{Min}) do satisfy the properties
(\ref{copM}).
We have already shown that
the functionals $f$ in (\ref{fin}) satisfy
(\ref{propf1}) ($\epsi (f^i_{~j} )=\de^i_j $
obviously also holds for this subset).

Consider now the last property (\ref{propM}).
We know that it holds for $\SOqrNt$.  Take  the free indices
$j$ and $k$  as {\small $\bu b$, $\bu\bu$} and {\small $\bu\ci$},
and apply the projection $P$ on both members of the equation.
It is an easy matter
to show that  only the $f$'s in (\ref{fin}) and the
$M$'s in (\ref{Min}) enter the sums: this
is due to the vanishing of some $P(M)$ and to some $f$'s.
We still have to prove that the $*$ product
in (\ref{propM})
can be computed via the coproduct $\D$ in
$\ISOqrN$.
Consider the projection of property
 (\ref{propM}), written symbolically as:
\eq
P [ M (f \otimes id) \D_{N+2} (a)]=P [(id\otimes f)
\D_{N+2} (a) M ]~. \label{PMfa}
\en
Now apply the definition (\ref{deffb}) and the first of
(\ref{co-iso}) to rewrite
(\ref{PMfa}) as
\eq
P(M)(\fb\otimes id)\D(P(a))=(id\otimes\fb)\D(P(a))P(M)~.
\en
This projected equation then  {\sl becomes}
 property (\ref{propM})  for the $\ISOqrN$
functionals
$f$ and adjoint elements $M$, with the
correct coproduct $\D$ of $\ISOqrN$.  Q.E.D.
\sk
To simplify notations, we write the composite indices
as follows:
\eq
{}_{\bu} {}^a \rightarrow {}^a,~{}_{\bu} {}^{\bu} \rightarrow
{}^{\bu},~{}_{\bu} {}^{\ci} \rightarrow {}^{\ci};~~~{}^{\bu} {}_a
\rightarrow
{}_a,~{}^{\bu} {}_{\bu} \rightarrow {}_{\bu},~{}^{\bu} {}_{\ci}
\rightarrow
{}_{\ci}
\en
Similarly we'll write $q_b$ instead of  $q_{ b \bu}  $.
\sk
\noi Using the general formula (\ref{omia}) we can
deduce the $\om, T$ commutations for $\ISOqrN$:
\eqa
& &\om^{b} \T{c}{d}= {q_{f} \over r} \Rinv{bf}{ed} \T{c}{f}
\om^{e}\label{VTcomm}\\
& &\om^{b} x^c={q_{b}\over r^2}  x^c \om^{b}
+ \lambda r^{\Nmezzi -1 } q_{d} C^{bd} \T{c}{d}
\om^{\ci}\label{Vxcomm}\\
& &\om^{b} u={r^2\over q_{b}} u~\om^{b}\\
& &\om^{b} v={q_{b}\over r^2} v~\om^{b}\\
& &\om^{\bu} \T{c}{d}= \T{c}{d}
\om^{\bu}\label{tauTcomm}\\
& &\om^{\bu} x^c={1\over r^2}  x^c \om^{\bu}
- \lambda {q_{b}\over r} \T{c}{b} \om^{b}   \label{tauxcomm}\\
& &\om^{\bu} u=r^2 u \om^{\bu}  \\
& &\om^{\bu} v=r^{-2} v \om^{\bu}  \\
& &\om^{\ci} \T{c}{d}= {q^2_{d} r^{-2} \T{c}{d} \om^\ci}
\label{omciTcomm}\\
& &\om^{\ci} x^c=x^c \om^{\ci}   \label{omcixcomm}\\
& &\om^{\ci} u= u \om^{\ci}  \\
& &\om^{\ci} v= v \om^{\ci}  \label{omvcomm}
\ena
\noi {\bf Note 4.1:} $u$ commutes with all
 $\om$ 's only if $r=q_{a}=1$ (cf. {\sl Note 2.2}). This means
that
$u=I$ is consistent with the differential calculus
on $ISO_{q_{ab};r=q_{a}=1}(N)$.
\sk
The 1-form $\tau \equiv \om^{\bu} \equiv \ome{\bu}{\bu}$
is bi-invariant, as one can check by using
 (\ref{DLin})-(\ref{DRin}).
Then an exterior derivative on $\ISOqrN$ can be defined
 as in eq. (\ref{defd1}), and satisfies the Leibniz rule.
The alternative expression $da=(\chi_i * a) \om^i$
(cf. (\ref{defd2}))
continues to hold, where
\eqa
& &\chi_b=\rinv \fmi{\bu}{b}   \nonumber \\
& &\chi_{\ci}=\rinv \fmi{\bu}{\ci} \nonumber \\
& &  \chi_{\bu}=\rinv[\fmi{\bu}{\bu}-\epsi]    \label{chiin}
\ena
are the left-invariant vectors dual to the left-invariant
1-forms $\om^b$, $\om^{\bu}$ and $\om^{\ci}$.
As a consequence of
(\ref{cofb}) their coproduct is given by
\eqa
& &\D (\chi_b)=\chi_{\bu}  \otimes \fmi{\bu}{b}+\chi_c \otimes
\fmi{c}{b} +
 \epsi \otimes \chi_b  \label{copchi1}\\
& &\D (\chi_{\bu}) = \chi_{\bu} \otimes \fmi{\bu}{\bu} + \epsi
\otimes
\chi_{\bu}
\label{copchi2}\\
& &\D (\chi_{\ci})=\chi_{\ci} \otimes \fmi{\ci}{\ci} + \chi_{\bu}
\otimes
 \fmi{\bu}{\ci} + \chi_{c}\otimes \fmi{c}{\ci}
+\epsi\otimes\chi_{\ci}
 \label{copchi3}
\ena
\sk
The exterior
derivative on the generators of $\ISOqrN$  reads:
\eqa
& & d\T{c}{d}=0 \label{trivial} \\
& & dx^c=-q_{b} r^{-1} \T{c}{b} \om^b - r^{-1}  x^c \om^{\bu}
\label{donx}\\
& & du = ru\om^{\bu} \\
& & dv = -r^{-1} v \om^{\bu} \label{donv}\\
& & dz=  - q_{b} r^{-1} y_b \om^b -r(1-r^N) u\om^{\ci} - r^{-1} z
\om^{\bu}
\label{donz}
\ena
\noi where we have included the exterior derivative on $z$ for
convenience.  Note that the calculus is trivial on the $\SOqrN$
subgroup of $\ISOqrN$, as is evident from (\ref{trivial}). Thus
effectively
we are discussing a bicovariant calculus on the orthogonal
$q$-plane generated by the coordinates $x^a$ and the
``dilatations" $u,v$.
\sk
Every element $\rho$ of $\Ga$ can be written as
$\rho=\sum_k a_k db_k$ for some $a_k,b_k$ belonging to
$\ISOqrN$. Indeed inverting  the  relations
(\ref{donx})-(\ref{donz})
yields:
\eqa
& & \om^a=-{r\over q_{a\bu}} \kappa (\T{a}{c})[ dx^c - x^c u dv]
\label{omdT1}\\
& & \om^{\bu} = -r u dv=r^{-1} v du \\
& & \om^{\ci}=-{{vdz + r^{-N}zdv + r^{-\Nmezzi} C_{ab} x^a
dx^b}\over
{r(1-r^N)}}
\label{omdT3}
\ena
Finally, the two properties (\ref{dbicov})
hold also for $\ISOqrN$, because of the bi-invariance
of $\tau=\om^{\bu}$.
Thus all the axioms for a bicovariant
first order differential calculus on $\ISOqrN$
are satisfied.
\sk
The exterior product of left-invariant one-forms
is defined as
\eq
\om^i \we \om^j\equiv \om^i \otimes \om^j - \L{ij}{kl}
\om^k \otimes \om^l \label{omcomgen}
\en
\noi where
\eq
\L{ij}{kl}=\f{i}{l} (\M{k}{j})
\en
This $\Lambda$ tensor can in fact
be obtained from
the one of $\SOqrNt$ by restricting its indices to the
subset  {\small $\bu b, \bu\bu, \bu\ci$}. This is true because
when  {\small $ i,l=  \bu b, \bu\bu, \bu\ci$ }  we have
$\f{i}{l} (Ker P)=0$ so that $\f{i}{l}$ is
well defined on $\ISOqrN$, and we can write
$\f{i}{l} (\M{k}{j})=\fb^i_{~l} [P(\M{k}{j})]$
(see discussion after
Theorem  4.1).
The non-vanishing components of $\Lambda$ read:
\[
\begin{array}{lll}
\LM{ad}{cb}={q_a \over q_{c}} r^{-1} \R{ad}{bc} &
 \LM{\bu\ci}{cb}= - {r^{-\Nmezzi-1}\over q_{c}} \lambda C_{bc} &
 \LM{\bu d}{c\bu}=r^{-2} \de^d_c \\
\LM{a\ci}{c \ci}=r^{-1} \lambda \de^a_c &
\LM{\ci d}{c \ci}=({r\over q_c})^2 \de^d_c &
\LM{a\bu}{\bu b}=\de^a_b \\
\LM{\bu d}{\bu b}=r^{-1} \lambda \de^d_b &
\LM{a\ci}{\ci b}=r^{-4} (q_a)^2 \de^a_b &
\LM{\bu\bu}{\bu\bu}=1 \\
\LM{a d}{\bu\ci}=-q_ar^{-\Nmezzi-1}\la C^{da} &
\LM{\bu\ci}{\bu\ci}=\la r^{-1}(1-r^{-N}) &
\LM{\ci \bu}{\bu \ci}=1 \\
\LM{\bu\ci}{\ci\bu}=r^{-4}&
\LM{\ci\ci}{\ci\ci}=1& \end{array}
\]
{}From (\ref{omcomgen}) it is not difficult to deduce the
commutations
between the $\om$'s:
\eqa
& & {1\over q_c} \PS{ab}{cd} ~ \om^d \we \om^c=0\label{PSom}\\
& &  \om^a \we \om^{\bu} = -r^2 \om^{\bu} \we \om^{a} \\
& & \om^a \we \om^{\ci} = -r^{-4} (q_a )^2 \om^{\ci} \we \om^{a} \\
& & \om^{\bu} \we \om^{\bu} = \om^{\ci} \we \om^{\ci} = 0 \\
& & \om^{\bu} \we \om^{\ci}=-r^{-4} \om^{\ci} \we \om^{\bu}
+{\lambda r^{-\Nmezzi -1} \over  q_a (1-r^{-N})}C_{ba} ~\om^a \we
\om^b
\label{omomcomm}
\ena
where $P_S$ is the $q$-symmetrizer given in Appendix B.
Notice that the dimension of the space of  2-forms generated by
$\om^a \we \om^b$ is
larger than in the commutative case since $P_S$ project into an
$N(N+1)/2-1$
(and not into an $N(N+1)/2$) dimensional space.
This is not surprising since the exterior algebra of
homogeneous orthogonal quantum groups is known to be larger than its
classical counterpart.
\sk
The exterior differential on $\Ga^{\we n}$
can be defined as in Section 4 (eq. (\ref{defdgen})),
and satisfies all the properties (\ref{propd2}), (\ref{propd4}).
As for $\SOqrNt$ the last two hold because of the
bi-invariance of $\tau \equiv \om^{\bu}$.
\sk
The Cartan-Maurer equations
\eq
d\om^i=\rinv (\tau \we \om^i+\om^i \we \tau)
\en
can be explicitly found after use of the commutations  (\ref{PSom})-
(\ref{omomcomm}):
\eqa
& & d\om^a=r^{-1} \om^a \we \om^{\bu} \\
& & d\om^{\bu}=0 \\
& & d\om^{\ci}=-r(1+r^2) \om^{\bu} \we \om^{\ci} +{r^3 \over
r^{\Nmezzi} -
r^{-\Nmezzi}}  {C_{ba} \over q_a} \om^a \we \om^b
\ena
The nonvanishing structure constants $\Cb$ (appearing in the $q$-Lie
algebra , see below), given by $\C{jk}{i}=\chi_k(\M{j}{i})$,  read:
\[
\begin{array}{lll}
\C{ab}{\ci}=- q_a^{-1}  r^{-\Nmezzi - 1} C_{ba} &
\C{a\bu}{c} = -r^{-1} \de^c_a &
\C{\bu b}{c}=r^{-1} \de^c_b \\
\C{\ci\bu}{\ci}=-r^{-3} (1+r^2)&
\C{\bu\ci}{\ci}=r^{-1} (1-r^{-N})& \end{array}
\]
\noi  These structure constants
 can be obtained from those of $\SOqrNt$ by
specializing indices, for essentially the same reason
as for the $\Lambda$ components.
\sk
Using the values of the $\Lambda$ and $\Cb$
tensors given above,
we can explicitly write the $q$-Lie algebra of
translations and dilatations on $\ISOqrN$  as:
\eqa
& &\chi_{\ci}\chi_{b}- q_b^2 r^{-4} ~
\chi_{b}\chi_{\ci}=0\label{qlieprim}\\
& &\chi_{c}\chi_{\bu}-r^{-2}\chi_{\bu}
\chi_{c}=-r^{-1}\chi_{c}\label{qlieMin}\\
& &\chi_{\ci}\chi_{\bu}-r^{-4}\chi_{\bu}
\chi_{\ci}={-(1+r^2)\over {r^{3}}}\chi_{\ci}\\
& & q_{a} \PA{ab}{cd}\chi_{b} \chi_{a}=0
\label{quattro}
\ena
A combination of  the above relations yields:
\eq
\chi_{\ci} + \lambda ~\chi_{\bu} \chi_{\ci} = \lambda~
{-q_a r^{-\Nmezzi}\over {r^{-2} + r^{-N}}}~
 \chi_{a} C^{ba} \chi_{b}  \label{chibuci}
\en
Notice the similar structure of eq.s (\ref{PRTT44}) and
(\ref{chibuci}).
\sk
The *-conjugation on the $\chi$ functionals and on the one-forms
$\om$ can be deduced from (\ref{conjchiAB}) [use
$(q_{f})^{-1} \Dcal{f}{b}$ $={\overline{q}_{b}}\Dcal{f}{b}$]
\eqa
& &\!\!\!\!\!\!(\chi_b)^*=-r^{-N}
\Dcal{f}{b}{1\over q_f}D^d_{~f }\chi_d=-r^{-N} {\overline q_b}
\Dcal{f}{b}D^d_{~f}\chi_d=-r^{-N} {\overline q_b}
D^f_{~b}\Dcal{d}{f}\chi_d\label{*46}\\
& &\!\!\!\!\!\!(\chi_{\bu})^*=- \chi_{\bu}\\
& &\!\!\!\!\!\!(\chi_{\ci})^*=-r^{-2N-2}\chi_{\ci}\label{*48}
\ena
whereas the conjugation on the $\om$ one-forms can be deduced
from (\ref{om*chi}) and (\ref{*46})- (\ref{*48}) or directly from
their expression in terms of $dx,du,dv$ differentials
(\ref{omdT1})-(\ref{omdT3}) remembering
that $(da)^*=d(a^*)$:
\eqa
& &(\om^a)^*={\overline q_a}^{-1} r^{N} (D^{-1})^a_{~b}\Dcal{b}{c}
={\overline q_a}^{-1} r^{N} \Dcal{a}{b} (D^{-1})^b_{~c}
\om^c \\
& &(\om^{\bu})^*=\om^{\bu} \\
& &(\om^{\ci})^*=r^{2N+2} \om^{\ci}~.
\ena
%%%%%%%%%%%%%%%%%%%%%%%%%%%
\sect{Calculus on the multiparametric orthogonal quantum plane:
coordinates, differentials and partial derivatives}
%%%%%%%%%%%%%%%%%%%%%%%%%%%
\subsection{$ISO_{q,r}(N)$ bicovariant calculus}
In this section we concentrate on the
orthogonal quantum plane
\eq
M\equiv Fun_{q,r}\left( {ISO(N) \over SO(N)}\right)~,
\en
\noi i.e. the $ISO_{q,r}(N)$ subalgebra
generated by the coordinates $x^a$ and the dilatations
$u,v$.

We study  the action of the exterior differential $d$ on
$M$  and the corresponding space $\Ga_M$ of $1$-forms. $\Ga_M$ is
the sub-bimodule of $\Ga$ formed by all the elements $adb$ or
$(da')b'$
where $a, b, a', b'$ are polynomials in $x^a,u$ and $v$
[of course $adb=d(ab)-(da)b$].

We will see that a generic element $\rho$ of $\Ga_M$ cannot be
generated,
as a left module, only by the differentials $dx, dv$, i.e.
it cannot be written as $\rho=a_i dx^i+ a dv$. We need also to
introduce
the  differential $dz$ (or equivalently $dL\equiv d(x^aC_{ab}x^b)$).
Thus the basis of differentials  is
given by $dx^a,dv,dz$ and corresponds to the
intrinsic basis of {\sl independent} one-forms $\om^a, \om^{\bu}$
and
$\om^{\ci}$.
Note that $du$ can be expressed in terms of $dv$
since $du=-u(dv)u= -r^2u^2dv= -r^{-2} (dv)u^2$ [ see (\ref{udv})
below].
\sk
{\bf Commutations}
\sk
The commutations between the coordinates $x^a , u$ and $v$ have been
given
in Section 2.
The commutations between coordinates and differentials are found
by expressing the differentials in terms of the one-forms $\om$ as
in
(\ref{donx})-(\ref{donz}), and using then the $x,u,v$ commutations
with the $\om$'s given
in (\ref{VTcomm})-(\ref{omvcomm}). The resulting $q$-commutations
between $x$ and $dx$ are found to be:
\eq
(r^{-2}P_S-P_A)(x\otimes dx)=(P_S+P_A)(dx\otimes x)
\en
where the projectors $P_S$ and $P_A$ are defined in (\ref{projBCD}),
and we have used the tensor notation
$A^{ab}_{~~cd} x^c dx^d \equiv A (x\otimes dx)$
etc. The remaining commutations are
given in formulae (\ref{xdu})--(\ref{zdzcomm}) in Table 1.

Let us consider the above formula, giving the $x,dx$
commutations. If we
multiply it by $P_0$ we find $0=0$. Thus from this equation we have
no
information on $P_0 ( x \otimes dx) $. Applying instead the
projectors
$P_S$ and $P_A$ yields
\eq
P_S (x \otimes dx)= r^2 P_S (dx \otimes x)~;~~~
P_A (x \otimes dx)= - P_A (dx \otimes x)  \label{PSPAxdx}
\en
which does not  allow  to express $x^a dx^b$ only in terms of linear
combinations of  $(dx)x$ since no linear combination of $P_S$ and
$P_A$
is invertible. The space of $1$-forms has therefore one more
dimension
than
his classical analogue because we are missing a condition involving
the one dimensional projector
$P_0^{ab}{}_{\!\!\!ef}=Q_N(r)C^{ab}C_{~ef}$,
see (\ref{projBCD}).

However,  if we consider the $1$-form  $dL\equiv d(x^eC_{ef}x^f)$
-- an exterior derivative of {\sl polynomials} in the basic elements
--
we can  write  the
commutations
between
the $x$ and $dx$ elements  as follows:
\eqa
dx\otimes x
&=&-(P_S+P_A+P_0)x\otimes dx +(P_S+P_A)d(x\otimes x)+P_0d(x\otimes
x)
\nonumber\\
&=&P_S dx\otimes x +P_A dx\otimes x - P_0x\otimes dx +P_0d(x\otimes
x)
\nonumber\\
&=&(r^{-2}P_S-P_A-P_0)x\otimes dx + P_0d(x\otimes x)
\label{xdxcomm2}
\ena
where we have used the Leibniz rule, the commutations
(\ref{PSPAxdx})
and  $P_S+P_A+P_0=I$. Equivalently we have
\eq
dx \otimes x=(r^{-2}P_S-P_A-P_0) x\otimes dx
- C {r^{\Nmezzi-2} (1-r^2)\over 1-r^N} (v dz + zdv) \label{xdxcomm3}
\en
\noi involving the $dv$ and $dz$ differentials.
\sk
The presence of $dz$ can also be explained within the general theory
by recalling
that $\Ga$ is a free right module [see paragraph following
(\ref{eta})].
A basis of right invariant $1$-forms is given by (\ref{eta}) and
we explicitly have:
\eqa
\eta^{a}&=&-r^{-1} dx^a\,u=-r^{-1} d\T{a}{\bu}\,\kappa(\T{\bu}{\bu})
\label{right0} \\
\eta^{\bu}&=&-r^{-1} dv\,u =-r^{-1}
d\T{\bu}{\bu}\,\kappa(\T{\bu}{\bu})
\label{right1} \\
\eta^{\ci} &=& {r^{\Nmezzi-1}\over{(1-r^N)(1+r^{N-2})}}[dx^eC_{ef}x^f
-
                   r^{N-2} x^eC_{ef}dx^f] u^2
\label{righteta}\\
        &=&{-r^{N-1}\over {r^N-1}}[dz\,u
+dy_b\,\kappa(x^b)+du\,\kappa(z)]
         = {-r^{N-1}\over {r^N-1}} d\T{\ci}{B}\,\kappa(\T{B}{\bu})
\ena
To derive the expression for $\eta^{\ci}$ use:
$y_b=-r^{-\Nmezzi}ux^eC_{ef}\T{f}{b};$
$dy_b\,\kappa(x^b)=r^{-\Nmezzi}du\,xxu +r^{-\Nmezzi}dx\,xu^2;$
$dz=d({-r^{-\Nmezzi}\over{1+r^{2-N}}}uxx)$
$={-r^{-\Nmezzi}\over{1+r^{2-N}}}(du\,xx+
dx\,xu+xdx\,u)$;
$\kappa(z)=\kappa(\T{\ci}{\bu})=r^{-N}z$; $uxx=r^2xxu$;
$udx\,x=dx\,xu$,
where $xx\equiv L\equiv x^eC_{ef}x^f$.

The  $1$-forms (\ref{right1})-(\ref{righteta})
in $\Ga$ do not contain any $\T{a}{b}$ element
and therefore belong to $\Ga_M$ as well; they are linearly
independent
and freely generate $\Ga_M$ as a right module because they freely
generate
the full $\Ga$ as a right module. The extra $1$-form $\eta^{\ci}$ (or
$dz$)
is therefore a natural consequence of the right module structure of
$\Ga$.
\sk
In summary: either
$dL$
or $dz$ or $\eta^{\ci}$ are necessary in order to close the
commutation
algebra
between
coordinates and differentials. Thus the commutations involving $z$
and
$dz$ appear in Table 1.
\sk
We have seen that $dv\,u$; $dx^a\,u$ and $\eta^{\ci}$ freely
generate
$\Ga_M$ as
a right module; recalling that $\Ga$ is also a free left module, we
have the :
\sk
\noi {\bf Proposition 5.1} The $M$-bimodule $\Ga_M$, as a left
module
(or as a right module), is freely generated by the differentials
$dx$, $dL$ (or $dz$)
and $dv$.
\noi{\sl Proof:}  to show that
$a_idx^i+adL+a_{\bu}dv=0$ $\Rightarrow a_i=0, a=0, a_{\bu}=0$
express $dx^i, dL, dv$ in terms of $\om^a,\om^{\ci}\om^{\bu}$,
see (\ref{donx})-(\ref{donz}), and recall that $\Ga$ is a free left
module.
\sk
\noi {\bf Note 5.1}
{}From (\ref{xdu}), (\ref{xdz}) and the commutations of $L$ with $x$
and
$u$
we have $x^cdL=dL\,x^c$,
$udL=dL\,u$ and $vdL=dL\,v$.
These relations and (\ref{xdxcomm}) show that inside $\Ga_M$ there
is
the
smaller bimodule generated by the differentials $dx^a$ and $dL$.
\sk
We now  examine the space of 2-forms.
By simply applying the exterior derivative $d$ to the relations
(\ref{xdxcomm})-(\ref{zdzcomm})  we deduce
the commutations between the differentials given in Table 1.
As with the $\om^a$'s in eq. (\ref{PSom}), the relations in
(\ref{dxdxcomm})
are not
sufficient to order the differentials $dx^a$.
\sk
{\bf $ISO_{q,r}(N)$ - coactions}
\sk
All the relations we have been deriving have
many symmetries properties because they are covariant
under  the actions on $M$ and $\Ga$ of the full  $ISO_{q,r}(N)$
$q$-group. In fact we have the following three $ISO_{q,r}(N)$
actions:

1) the coproduct of $ISO_{q,r}(N)$ can be seen as a left-coaction
$\D :~M\rightarrow ISO_{q,r}(N)\otimes M$:
\eq
\D x^a= \T{a}{b}\otimes x^b + x^a\otimes v\label{1)}
\en

2) the left coaction
$\DL:\Ga\rightarrow ISO_{q,r}(N)\otimes \Ga$,
when restricted to $\Ga_M$ gives
\eq
\left.\DL\right|_{{}_{\Ga_M}}:\Ga_M\rightarrow ISO_{q,r}(N)\otimes
\Ga_M\label{2)}
\en
and defines a left coaction of $ISO_{q,r}(N)$ on $\Ga_M$ compatible
with
the bimodule structure of $\Ga_M$ and the exterior differential:
$\left.\DL\right|_{{}_{\Ga_{\!M}}}(adb)
=\D(a)(id\otimes d)\D(b)$.

3) the right coaction $\DR:\Ga\rightarrow \Ga\otimes ISO_{q,r}(N)$
does not become a right coaction of $ISO_{q,r}(N)$ on
$\Ga_M$; however we have
\eq
\left.\DR\right|_{{}_{\Ga_{\!M}}}:\Ga_M\rightarrow \Ga\otimes
M\subset
\Ga\otimes ISO_{q,r}(N)\label{3)}
\en
this map is obviously well defined and satisfies
$\left.\DR\right|_{{}_{\Ga_{\!M}}} (a db)=
\D(a)(d\otimes id)\D(b)$
%$(id\otimes\D)\left.\DR\right|_{{}_{\Ga_{\!M}}}(\rho)=
%(\DR\otimes id)\left.\DR\right|_{{}_{\Ga_{\!M}}}(\rho)$;
%$(id\otimes \epsi)\left.\DR\right|_{{}_{\Ga_{\!M}}}(\rho)=\rho$
%and (\ref{dbicov}) $\forall \rho\in \Ga_M$ and
$\forall a, b \in M$ since $M\subset ISO_{q,r}(N)$.

We call this calculus $ISO_{q,r}(N)$-bicovariant
because $\left.\DL\right|_{{}_{\Ga_{\!M}}}$ and
$\left.\DR\right|_{{}_{\Ga_{\!M}}}$ are compatible with the
bimodule structure
of $\Ga_M$ and with the exterior differential.
\sk
{\bf Partial derivatives}
\sk
The tangent vectors $\chi$ in (\ref{chiin})
and the corresponding  vector fields $\chi *$
 have ``flat" indices. To compare $\chi *$
with partial derivative operators with ``curved" indices,
we need to define the operators $\stackrel{\leftarrow}{\chit} :$
\eqa
& &\stackrel{\leftarrow}{\chit_c}\!(a) \equiv - {r\over q_{b}}
(\chi_b * a) \kappa (\T{b}{c})-{r^{-\Nmezzi} \over r(1-r^N)}~
(\chi_{\ci} *a) C_{bc}x^b, \label{partialcurved1}\\
& &\parleft_{\bu} (a)\equiv  {r\over q_{b}}
(\chi_b * a) \kappa (\T{b}{c}) x^c u-r(\chi_{\bu} *a)u-
{r^{-N} \over r(1-r^N)}~
(\chi_{\ci}*a) z \label{partialcurved2}\\
& &\parleft_{\ci} (a) \equiv - {1 \over r(1-r^N)} (\chi_{\ci} *a) v
\label{partialcurved3}
\ena
\noi so that
\eq
d\,a= \stackrel{\leftarrow}{\chit_c}\!(a) ~dx^c + \parleft_{\bu} (a)
dv
+ \parleft_{\ci} (a) dz \equiv \parleft_{C} (a) dx^C
\label{dacoord1}
\en
\noi [{\small $C=(\ci,c,\bu)$}, $dx^C = (dz,dx^c,dv)$] which is
equation
$da=(\chi_c * a) \om^c + (\chi_{\bu} *a) \om^{\bu}
+ (\chi_{\ci} * a) \om^{\ci} $  in ``curved" indices.
The action of $\stackrel{\leftarrow}{\chit_C}$ on the coordinates
$x^C = (z,x^c,v)$ is given by
\eq
\stackrel{\leftarrow}{\chit_C}\!(x^A)=\de^A_C I~,
\en
{}From the Leibniz rule $d(ab)=(da)b+a(db)$,
 using (\ref{dacoord1}) and the
fact that  $dx^C = (dz,dx^c,dv)$  is a basis for 1-forms, we
find for example:
\eq
\parleft_c (ax^b)=a\delta^b_c + \parleft_d (a)
(r^{-2}P_S-P_A-P_0)^{db}_{~~ec}
x^e-(1-r^2) \parleft_{\bu}(a) \de^b_c v \label{Leibpart}
\en
The tangent vector fields $\chi_c *$ of this paper and the partial
derivatives $\parleft$
are derivative operators that act  ``from the right to the left"
as can be seen from the deformed
Leibniz rule (\ref{Leibpart}). This explains the
inverted arrow on $\parleft$.

Eq. (\ref{Leibpart}) gives the $\parleft_c, x^b$ commutations. The
rest of the $\parleft_C, x^B$ commutations reads:
\eqa
& &\parleft_{\bu} (ax^b)=q_b^{-1} \parleft_{\bu} (a) x^b -
C^{cb} {r^{\Nmezzi-2}(1-r^2) \over (1-r^N)}~ \parleft_c (a) z \\
& &\parleft_{\ci} (ax^b)=q_b \parleft_{\ci} (a) x^b -
C^{cb} {r^{\Nmezzi-2}(1-r^2) \over (1-r^N)}~ \parleft_c (a) v
\ena
\eqa
& &\parleft_c(av)=r^{-2} q_c \parleft_c(a) v\\
& &\parleft_{\bu} (av)= r^{-2} \parleft_{\bu} (a) v+a\\
& &\parleft_{\ci} (av)=\parleft_{\ci} (a) v
\ena
\eqa
& &\parleft_c(az)= q_c^{-1} \parleft_c(a) z\\
& &\parleft_{\bu} (az)= r^{-2} \parleft_{\bu} (a) z\\
& &\parleft_{\ci} (az)=r^2 \parleft_{\ci} (a) z +a+(r^{-2}-1)
\parleft_{\bu} (a) v
\ena
We can also define derivative operators acting from the left
to the right, as in ref.s \cite{qplane},
using the antipode $\kappa$ which is antimultiplicative
[one can also use (\ref{defd2eta})].
For a generic quantum group the vectors ${-\kpm}(\chi_i) \equiv
-\chi_i \ci
\km$
act from the left and we also have
\eq
d\,a=(\chi_i*a)\omega^i=\omega^i({-\kpm}(\chi_i)*a)
\en
as is seen  from $\kp (\chi_i)=-\chi_j \kp (\f{j}{i})$ and $\kpm
(\f{k}{j})
\f{i}{k}
= \de^i_j$ [third line of (\ref{copf})].

We then define the partial derivatives $\partial_C$ so that
\eq
d\,a=dx^C\;\chit_C(a)~.
\en
Again the action  of $\chit_C$ on the coordinates is
\eq
\chit_C (x^A)=\de^A_C I~,
\en
The $\partial_C,x^B$ commutations are given in Table 1.

\subsection{$SO_{q,r}(N)$ bicovariant calculus}

{\bf Commutations}
\sk
Since the $\PA{ab}{cd} x^c x^d=0 $ commutation relations allow for
an
ordering of the coordinates (moreover the Poincar\'e series of the
polynomials on the quantum orthogonal plane is the same as the
classical
one), it is tempting to impose extra conditions on the differential
algebra
of the $q$-Minkowski plane, so that the space of $1$-forms has the
same
dimension as in the classical case. We require that the
commutation relations between $x$ and $dx$ close on the algebra
generated
by $x$ and $dx$:
\eq
dx^a x^b= \alpha^{ab}{}_{\!\!ef}x^edx^f \label{dxxalfa}
\en
where $\alpha$ is an unknown matrix whose entries are complex
numbers.
Any such matrix can be expanded as $\al=aP_S + bP_A + cP_0$ with
$a,b,c=const$.
{} From (\ref{PSPAxdx}) we have $\alpha=r^{-2}P_S-P_A+cP_0$;
 therefore condition (\ref{dxxalfa})
is equivalent to
\eq
P_0(dx\otimes x)=
 c P_0(x\otimes dx)\label{P0dxx}
\en
and supplements eq.s  (\ref{PSPAxdx}). Taking its exterior
derivative
leads to a supplementary condition on the $dx,dx$ products (for $c
\not= -1$):
\eq
P_0 (dx \we dx)=0 \label{P0dxdx}~.
\en
{}From  (\ref{dxdxcomm}) and (\ref{P0dxdx}) it follows that
$dx \we dx = (P_S + P_A + P_0) (dx \we dx )= P_A (dx \we dx)$, or
[see the definition of $P_A$ in (\ref{projBCD})] :
\eq
dx\we dx=-r\Rha ~dx\we dx\label{dxdxcomm2}~.
\en
\noi which allows the ordering of $dx,dx$ products.

Using (\ref{xdxcomm2}), (\ref{P0dxx}) and (\ref{RprojBCD}), we find
\eqa
dx\otimes x
& =& (r^{-2}P_S-P_A) (x\otimes dx) + P_0 (dx \otimes x)\nonumber\\
& =& (r^{-2}P_S-P_A) (x\otimes dx) + c P_0 (x \otimes dx)\nonumber\\
&=& (r^{-2}P_S-P_A+r^{N-2}P_0) (x\otimes dx) + (c-r^{N-2}) P_0 (x
\otimes
dx)\nonumber\\
&=& r^{-1}{\hat{R}}^{-1}(x\otimes dx)
+ (c-r^{N-2})P_0(x\otimes dx)
\label{xdxP0}~.
\ena
The consistency of the
commutation relations (\ref{dxdxcomm2}) and (\ref{xdxP0}) with the
associativity condition on the triple $dx^i\,dx^j\,x^k$ fixes
 $c=r^{N-2}$ i.e.:
\eq
P_0(dx\otimes x)=
r^{N-2} P_0(x\otimes dx)\label{P0dxx2}~;
\en
the  $x,dx$ commutations (\ref{xdxP0}) then become:
\eq
x\otimes dx = r \Rha (dx \otimes x) \label{xdxcomm4}
\en
\noi and reproduce (in the uniparametric case) the
known $x,dx$ commutations of the quantum orthogonal plane
\cite{Fiore}.
\sk
{\bf Coactions}
\sk
This calculus is no more covariant under the $\ISOqrN$ action,
\eq
x^a \longrightarrow \T{a}{b} \otimes x^b~+x^a\otimes v~, u
\longrightarrow u
\otimes u,~v \longrightarrow  v \otimes v \label{ISOcov}
\en
but we are left with covariance under the $SO_{q,r}(N)$ action
\eq
x^a \longrightarrow \T{a}{b} \otimes x^b~. \label{SOcov}
\en
In other words, $\delta_L : \Ga_M'\rightarrow SO_{q,r}\otimes
\Ga_M'$
defined by $\delta_L(adb)=\delta(c)(id\otimes d)\delta(b)$
with $\delta(x^a)=\T{a}{b}\otimes x^b$ is a left coaction of
$SO_{q,r}(N)$ on the bimodule $\Ga_M'$ where $\Ga_M'$ is $\Ga_M$
with the extra condition (\ref{P0dxx}) [cf. (\ref{2)})].
Similarly, the map $\delta_R(adb)=\delta(a)(d\otimes id)\delta(b)$
is well defined [cf. (\ref{3)})].

Left covariance under (\ref{ISOcov}) is broken only by
(\ref{P0dxx}).
Indeed, while relations (\ref{PSPAxdx}) are left and right
$ISO_{q,r}(N)$-covariant,
the extra condition (\ref{P0dxx}) is {\sl not}
left $ISO_{q,r}(N)$-covariant :
$\DL[P_0(dx\otimes x)- c P_0(x\otimes dx)]\not=0, \forall c$. It is
right $ISO_{q,r}(N)$-covariant,
$\DR[P_0(dx\otimes x)- c P_0(x\otimes dx)]=0$,
only for $c=r^{N-2}$, as can be seen
using
$\T{b}{d}dx^a=d(\T{b}{d}x^a)={r\over {q_{d}}}R^{ab}_{~ef}dx^e\,
\T{f}{d}$
and (\ref{CR}).
Therefore the choice $c=r^{N-2}$ preserves  the right
coaction
$\DR$.

\sk
\noi {\bf Note 5.2}  We can reformulate the quotient
procedure $\Ga_M\rightarrow \Ga_M'$
in a more abstact setting by considering that $\Ga_M$ is a
subbimodule of
the {\sl bicovariant} bimodule $\Ga$.
In  (\ref{righteta}) we have expressed the
$x^eC_{ef}dx^f\leftrightarrow
dx^eC_{ef}x^f$ commutation via the right invariant $1$-form
$\eta^{\ci}$.
A condition on $\Ga$ (and therefore on $\Ga_M$) that preserves
{\sl right $ISO_{q,r}(N)$ covariance}, i.e. compatible with $\DR$
[as
given
in (\ref{DRin})], is: $\eta^{\ci}$ linearly dependent
from the remaining right invariant $1$-forms $dv\,u$ and $dx^a\,u$.
It is
easily seen that since $\eta^{\ci}$ is quadratic in the basis
elements
$x^a$ the
only possible linear condition is $\eta^{\ci}=0$, and this gives
exactly
(\ref{P0dxx2}).
The $M$-bimodule $\Gamma_M'$ is therefore
generated by the differentials $dx^b$ and
$dv$. Since left $ISO_{q,r}(N)$ covariance is broken  (whereas  right
$ISO_{q,r}(N)$ covariance is preserved), the relation between the
left
invariant $1$-forms is {\sl nonlinear}. Explicitly we have
\eq
\om^{\ci}=-{q_a \over r^2} vy_a ~\om^a+{r^{\Nmezzi-2} \over r^N +
r^2}
C_{ab} x^a x^b \om^{\bu}
\en
\noi [express $dz$ in terms of $dx^i$,$dv$ in (\ref{omdT3}) and use
the expansion of $dx^b$ and $dv$ on $\om^a$ and $\om^{\bu}$ as given
in (\ref{donx}),(\ref{donv})].
\sk
{\bf Partial derivatives}
\sk
The relevant $\parleft$ operators
reduce to $\parleft_c, \parleft_{\bu}$, and the $\parleft_C, x^b$,
commutations become:
\eqa
& &\parleft_c (ax^b)=a\delta^b_c + \parleft_d (a) r^{-1}
(\Rha^{-1})^{db}_{~~ec}
x^e-(1-r^2) \parleft_{\bu} \de^b_c v \label{Leibpart2}\\
& &\parleft_{\bu} (ax^b)=q_b^{-1} \parleft_{\bu} (a) x^b
\ena
\noi while the $\parleft_C, v$ commutations are unchanged.
Note the dilatation operator $\parleft_{\bu}$ appearing on the
right-hand side of (\ref{Leibpart}) or (\ref{Leibpart2}).
The $\partial_C,x^B$ commutations with  {\small
$C=(c,\bullet),\,B=(b,\bullet)$}
are collected in Table 2.

{}From $d^2(a)=0=d(\parleft_C (a) dx^C)=\parleft_B (\parleft_C (a))
dx^B \we dx^C$ and the $q$-commutations of the differentials
(\ref{dxdxcomm})-(\ref{dzdzcomm})
one finds the commutations between the partial
derivatives:
\eqa
& & (P_A)^{ab}_{~~cd} \parleft_a \parleft_b=0  \label{parleftcomm} \\
& & \parleft_b \parleft_{\bu} - { q_b\over r^2 } \parleft_{\bu}
\parleft_{b}=0 \\
& & \parleft_b \parleft_{\ci}-q_b \parleft_{\ci} \parleft_b=0 \\
& & \parleft_{\bu} \parleft_{\ci}-\parleft_{\ci} \parleft_{\bu}=0
\ena
Similarly the $\partial_A, \partial_B$ commutations are given in
Table 2.
\sk
We now give an explicit relation between the
$\partial_c$, $\partial_{\bu}$ and the $q$-Lie algebra generators
$\chi_{{}_{c}}$, $\chi_{{}_{\bu}}$
(a similar expression holds also for the $\parleft$ derivatives).
{}Recalling (\ref{defd2eta}) we have:
\eq
d\,a=-\eta^C(a*\kp(\chi_{{}_C}))\label{etakchi}
\en
where here {\small $C=(c,\bu)$} because we have set $\eta^{\ci}=0$.
Putting together (\ref{etakchi}) and (\ref{right0}),(\ref{right1}),
the
relations $d\,a=dx^C\;\chit_C(a)$
give, $\forall a\in \ISO$  :
\eq
\chit_c(a)=r^{-1}u(a*\kp(\chi_c))~~~,~~~~~~
\chit_{\bu}(a)=r^{-1}u(a*\kp(\chi_{\bu}))
\label{pechi}~.
\en
The commutations between the partial derivatives
given in Table 2
were obtained from $d^2=0$, but
can be also derived
via (\ref{pechi}) and the $q$-Lie algebra
(\ref{qlieprim})-(\ref{quattro}).

Similarly we can introduce the right invariant
vector fields
\eq
h_C \equiv
h_{\kp(\chi_C)} \equiv
[\kp(\chi_{{}_C})\otimes id] \D
\en
and use their Leibniz rule
[which follows from
$\D(\kp(\chi_{{}_C}))=\kp(\chi_{{}_C})\otimes\epsi +
\kp(\f{D}{C})\otimes \kp(\chi_{{}_D})\,$]:
\eq
h_C(ab)=h_C(a)\,b+\kp(\f{D}{C})(a_1)\,\,a_2h_D(b)
\en
to rederive the $\partial, x,u$ commutations.
For example we have $h_ax^b=rv\delta^b_a+(r/q_b) R^{eb}_{~ac}x^c
h_e$ $+r\lambda h_{\bu}$ that together
with $\partial_C=r^{-1}uh_C$
(cf. eq.s (\ref{pechi}))
gives
\[\chit_ax^b=\delta^i_j[I+(r^2-1)v\partial_{\bu}]
+rR^{eb}_{~ac}x^c\chit_e~.
\]

\sk
{\bf Conjugation}
\sk
The commutations in Table 2 are consistent under the conjugation
(already defined for $x^a$ and $dx^a$)
\eq
(x^a)^*=\Dcal{a}{b}x^b,~(dx^a)^*=\Dcal{a}{b}dx^b,~(\partial_a)^*=
-r^Nd^{-1}_b \Dcal{b}{a}  \partial_b
\label{dstar1}
\en
\eq
v^*=v,~ (dv)^*=dv,~(\partial_{\bu})^*=u- \partial_{\bu}
 \label{dstar2}
\en
\noi where the entries $d_a$ have been defined after eq.
(\ref{Lambda}).
This can be proved directly by taking the *-conjugates of
the relations in Table 2, and by using the identity (\ref{Rprop2})
and:
\eq
{\bar C}=C^T;~ [Q_N(r)]^*=Q_N(r);
\en
\eq
\begin{array}{ll}
{}~d^c d^{-1}_h R^{cg}_{~~ha} (R^{-1})^{ea}_{~~cd}=
\de^e_h \de^g_d; &~ R^{ab}_{~~cd} d^a d^b=R^{ab}_{~~cd} d^c d^d
\nonumber\\
d^c R^{cg}_{~~hc}=r^{N-1} \de^g_h; &~
R^{ab}_{~~cd} d^{-1}_a d^{-1}_b=R^{ab}_{~~cd} d^{-1}_c d^{-1}_d
\nonumber
\end{array}
\en
\eq
\overline{q_a}={1 \over q_a} \mbox{ for $a
\not=n,n+1$},~~\overline{q_n}=
{1\over q_{n+1}}
\en

We now {\sl derive}
the conjugation on the partial derivatives from the differential
calculus on
$\ISO$. This is achieved by studying the conjugation on the right
invariant
vector fields $h$.
\sk
For a generic Hopf algebra, with tangent vectors $\chi_i$, we
deduce
the conjugation on $h$ from the commutation relations between
$h$ and a generic element of the Hopf algebra:
\eq
h_jb=h_j(b)+\kp(\f{s}{j})(b_1)\,
\,b_2h_s~={{}} \kp(\chi_j)(b_1)\,\, b_2+{{}}\kp(\f{s}{j})(b_1)\,
\,b_2h_s~
\en
We multiply this expression by ${{}} {\kp}^2(\f{j}{i})(b_0)\,$
[where
we have used the notation $(id\otimes \D)\D(b)=b_0\otimes b_1\otimes
b_2$]
to obtain
\eq
{{}}{\kp}^2(\f{j}{i})(b_1)\,\,h_jb_2+
{{}}{\kp}^2(\chi_i)(b_1)\, b_2=bh_i
\en
Now, using ${{}}\psi(b)\,=\overline{{{}}[\kp(\psi)]^*(b^*)\,}$
and then applying $*$ we obtain (here $a=b^*$)
\eq
h_j^*a={{}} [{\kp}^3(\chi_j)]^*(a_1)\,\,
a_2+{{}}[{\kp}^3(\f{s}{j})]^*(a_1)\,
\,a_2h_s~.
\en
This last relation implies
\eq
h_i^*\equiv [h_{\kp(\chi_i)}]^*
=h_{[\kp^3(\chi_i)]^*}\label{genfor}
\en
Notice that $*\sma{$\circ$}\kp^2$ is a  well defined conjugation
since
$(*\sma{$\circ$}\kp^2)^2=id$ .
\sk
We now apply  formula (\ref{genfor}), valid for a generic Hopf
algebra,
to the  $*$-conjugation
and the right invariant vector fields  of this section; we have:
\eqa
& &[h_{\kp(\chi_a)}]^*=-r^N\overline{q}_ad^{-1}_a\Dcal{b}{a}
h_{\kp(\chi_b)}\\
& &[h_{\kp(\chi_{\bu})}]^*=-
h_{\kp(\chi_{\bu})}~.
\ena
{} From these last relations and $\partial_C=r^{-1} u h_C$ we
finally deduce
$
(\partial_a)^*=
-d^{-1}_b \Dcal{b}{a} r^{N} \partial_b$ and
$(\partial_a)^*=(\partial_{\bu})^*=u- \partial_{\bu}~$ as in  
(\ref{dstar1}),
(\ref{dstar2}).

\subsection{The reduced $SO_{q,r}(N)$-bicovariant  algebra generated
by $x^a,
dx^a, \partial_a$ }
\sk
Note that the algebra in Table 2 actually contains a subalgebra
generated only by $x^a,dx^a, \partial_a$:  indeed $\partial_{\bu}$
vanishes when acting on monomials containing only the coordinates
$x^b$,
as can be seen from (\ref{Leibpart5}).
This calculus is $\ISO$-right covariant because it can also be
obtained
imposing the conditions $\eta^{\bu}=0$ and $\chi_{\bu}=0$ that are
compatible
with the right coaction $\DR$ and the bimodule structure given by the
$\f{i}{j}$
functionals.
\sk
Table 3 contains the multiparametric orthogonal quantum plane
algebra
of coordinates, differentials and partial derivatives, together with
a consistent conjugation. We emphasize here
that this conjugation {\sl does not } require an additional scaling
operator
as in ref. \cite{Lorek}. Thus the algebra in Table 3 can be taken
as
starting point for a deformed Heisenberg algebra (i.e. a deformed
phase-space) .
\sk
{\bf Real coordinates and hermitean momenta}
\sk
We note that for the real form $ISO_{q,r}(n+1,n-1)$, the
transformation
\eqa
& &X^a={1\over \sqrt{2}}(x^a+x^{a'})\,,~~~\sma{$a \leq
n$}\label{54ii}\\
& &X^{n+1}={i\over \sqrt{2}}(x^n - x^{n+1})\\
& &X^a={1\over \sqrt{2}}(x^a-x^{a'})\,,~~~\sma{$a >
n+1$}\label{56ii}
\ena
\noi defines real coordinates $X^a$.
 On this basis the metric becomes $C'=(M^{-1})^T C M^{-1}$
(where $M$ is the transformation matrix $X=Mx$):
\eq C'=
{\small
{1\over 2}  \left (\matrix{
r^{\Nmezzi-1}+r^{-\Nmezzi+1}& 0 &{}&{}&0&
-(r^{\Nmezzi-1}-r^{-\Nmezzi+1})\cr
0 & r^{\Nmezzi-2}+r^{-\Nmezzi+2}&{}&{}&
-(r^{\Nmezzi-2}-r^{-\Nmezzi+2})&0 \cr
 \vdots & \vdots  &{}&{}&\vdots & \vdots\cr
0&0&2&0&0&0\cr
0&0&0&2&0&0\cr
\vdots & \vdots &{}&{} & \vdots & \vdots\cr
0 & r^{\Nmezzi-2}-r^{-\Nmezzi+2} &{}&{}&
-(r^{\Nmezzi-2}+r^{-\Nmezzi+2})&0 \cr
r^{\Nmezzi-1}-r^{-\Nmezzi+1}& 0 &{}&{}&0&
-(r^{\Nmezzi-1}+r^{-\Nmezzi+1})\cr
}\right ) }
\label{Cprimo=}
\en
and reduces for $r\rightarrow 1$ to the usual $SO(n+1,n-1)$
diagonal metric with $n+1$  plus signs and $n-1$ minus signs.
Notice that the diagonal elements of $C'$ are real while the off
diagonal ones
are imaginary; moreover $C'$ is hermitian (and can therefore be
diagonalized via a unitary matrix).
\sk
As for the coordinates $X$, it is possible to define antihermitian
$\chi$ and
$\partial$, and real
$\om$ and $dx$.
To  define hermitian momenta
we first notice that the partial derivatives
\eq
\tilde\partial_a\equiv r^{N\over 2} d_a^{-{1\over 2}}\partial_a
\en
behave,  under the hermitian conjugation $*$, similarly to the
coordinates $x^a$ (see (\ref{dstar1}), (\ref{dstar2})):
\eqa
& &(\tilde{\partial}_a)^* =-\tilde{\partial}_a~~~~~\sma{$a\not=
n,\,n+1$}\\
%~~~\mbox{if }~\sma{$a\not= n,n+1$}\\
& &(\tilde{\partial}_n)^* =-\tilde{\partial}_{n+1}~.
\ena
As in (\ref{54ii})--(\ref{56ii}) we then define:
\eqa
& &P_a={-i\hbar\over \sqrt{2}}(\tilde\partial_a+\tilde\partial_{a'})
\,,~~~\sma{$a \leq n$}\\
& &P_{n+1}={-\hbar\over
\sqrt{2}}(\tilde\partial_n-\tilde\partial_{n+1})\\
& &P_a={-i\hbar\over \sqrt{2}}(\tilde\partial_a-\tilde\partial_{a'})
\,,~~~\sma{$a > n+1$}
\ena
It is easy to see that the $P_a$ are hermitian: ${P_a}{}^*=P_a$, and
that
in the classical limit are the momenta conjugated to the coordinates
$X^a$:
$P_a(X^b)=-i\hbar\delta_a^b$. In the $r\not=1$ case we explicitly
have
(use $d_{a'}=d_a^{-1}\,,~ d_n=d_{n+1}=1$):
\[
\begin{array}{l}
P_a(X^a)=-{1\over 2} i \hbar r^{N\over 2} (d_a^{1\over
2}+d_a^{-{1\over
2}})\\[1em]
P_a(X^{a'})=- {1 \over 2} i \hbar \epsilon_a r^{N\over 2}
(d_a^{1\over
2}-d_a^{-{1\over 2}})
{}~~\sma{ where $\epsilon_a=1$ if $a < n$ and
$\epsilon_a=-1$ if $a > n+1$}
\end{array}
\]
while the other entries of the  $P_a(X^b)$ matrix are zero.
\sk
By defining the transformation matrix $N_{a}^{~b}$ as:
\eq
P_a \equiv -i\hbar N_a^{~b} \partial_b
\en
\noi  we find the deformed canonical commutation relations:
\eq
P_a X^b -r S^{bc}_{~~ad} X^d P_c = - i \hbar E_a^b I \label{PXcomm}
\en
\noi where
\eq
S^{bc}_{~~ad}=N_a^{~e} M^b_{~f} \Rhat{fh}{eg} (M^{-1})^g_{~d}
(N^{-1})_h^{~c},
{}~~~E_a^b \equiv {i\over \hbar} P_a(X^b) = N_a^{~c} M^b_{~c}
\en
Similarly one finds all the remaining commutations of the
$P$, $X$ and $dX$ algebra.
Notice that no unitary operator appears on the right-hand
side of  (\ref{PXcomm}). Our conjugation is consistent with
(\ref{PXcomm})
without the need of the extra operator of  ref.  \cite{Lorek} .
\sk
For $n=2$ the results of this section immediately yield the
bicovariant calculus on the quantum Minkowski
space, i.e. on the
multiparametric orthogonal quantum plane
$Fun_{q,r}(ISO(3,1)/SO(3,1))$. The relevant formulas
are collected in Appendix A.
 \vskip 1cm
\noi {\large \bf Table 1: the $ISO_{q,r}(N)$-bicovariant
$x^A,\partial_A,dx^A$
algebra}
\sk
\eqa
& &\PA{ab}{cd} x^c x^d=0 \\
& &x^b v=q_{b} v x^b~;~~~x^b u=q^{-1}_{b} u x^b \\
& &z=-{1\over {(r^{-\Nmezzi}+ r^{\Nmezzi-2})}} x^b C_{ba} x^a u \\
& &zv=r^2vz~;~~~zu=r^{-2}uz\\
& &q_a x^a z = z x^a
\ena
\eqa
& &
(x \otimes dx)=  (r^2 P_S-P_A-P_0 )(dx \otimes x) +P_0 d(x\otimes x)
 \label{xdxcomm}\\
& & x^c du= {1\over q_c} (du) x^c - {\la \over r} (dx^c) u ;~~
x^c dv= q_c (dv)x^c +\la r (dx^c) v \label{xdu}\\
& & x^c dz={1\over q_c} (dz) x^c\label{xdz}\\
& & \nonumber\\
& & u dx^c = {q_c \over r^2} (dx^c) u\label{udx}\\
& & u du=r^{-2} (du) u;~~u dv = r^{-2} (dv) u\label{udv}\\
& & u dz=(dz) u\label{udz}\\
& &\nonumber\\
& & v dx^c = {r^2 \over q_c } (dx^c) v\label{vdx}\\
& & v du=r^{2} (du) v;~~v dv = r^{2} (dv) v\label{vdu}\\
& & v dz=(dz) v\\
& & \nonumber\\
& & z dx^c = q_c (dx^c) z \label{zdx}\\
& & z du =  r^{-2} (du) z + (r^{-2} - 1) (dz)u
\label{zdu}\\
& & z dv = r^2 (dv) z  + (r^2-1) (dz)v \label{zdv}\\
& & z dz=r^{-2} (dz) z \label{zdzcomm}
\ena
\eqa
& & P_S (dx \we dx) =0 \label{dxdxcomm}\\
& & dx^c \we du=-{r^2\over q_c} du \we dx^c;~~dx^c \we dv=-{q_c\over
r^2} dv
\we dx^c\\
& & dx^c \we dz=-{1\over q_c} dz \we dx^c \\
& & du \we du=dv \we dv = 0\\
& & du \we dv = - r^{-2} dv \we du =0\\
& & dz \we du = -du \we dz;~~dz \we dv = - dv \we dz\\
& & dz \we dz = 0 \label{dzdzcomm}
\ena
\eqa
& &\partial_c x^b=\delta^b_c I +  (r^{2}P_S-P_A-P_0)^{be}_{~~cd}
x^d \partial_e -(1-r^2) \de^b_c v \partial_{\bu} \label{Leibpart3}\\
& &\partial_{\bu} x^b=q_b x^b \partial_{\bu}  -
C^{bc} {r^{\Nmezzi}(1-r^2) \over (1-r^N)}~ z \partial_c  \\
& &\partial_{\ci} x^b=q_b^{-1} x^b \partial_{\ci}-
C^{bc} {r^{\Nmezzi}(1-r^2) \over (1-r^N)}~v \partial_c  \\
& &\partial_c v=r^{2} q_c^{-1} v \partial_c\\
& &\partial_{\bu} v= r^{2} v \partial_{\bu} +I\\
& &\partial_{\ci} v = v \partial_{\ci}  \\
& &\partial_c z= q_c z \partial_c\\
& &\partial_{\bu} z= r^{2} z \partial_{\bu}\\
& &\partial_{\ci} z=r^{-2} z  \partial_{\ci}  +I+(r^{2}-1)
v \partial_{\bu}
\ena
\sk\sk
\noi {\large \bf Table 2: the $SO_{q,r}(N)$ - bicovariant
$x^a,v,\partial_a,\partial_{\bu},dx^a,dv$ algebra}
\eqa
& &\PA{ab}{cd} x^c x^d=0 \\
& &x^b v=q_{b} v x^b
\ena
\eqa
& & x \otimes dx=  r\Rha (dx \otimes x)  \\
& & x^c dv= q_c (dv)x^c +\la r (dx^c) v \\
& & v dx^c = {r^2 \over q_c } (dx^c) v
\ena
\eqa
& & dx \we dx =-r {\hat R} dx \we dx \\
& & dx^c \we dv=-{q_c\over  r^2} dv \we dx^c\\
& & dv \we dv = 0
\ena
\eqa
& &\partial_c x^b= r \Rha^{be}_{~~cd} x^d \partial_e+
\delta^b_c [ I +(r^2-1) v \partial_{\bu}] \label{Leibpart4}\\
& &\partial_{\bu} x^b=q_b x^b \partial_{\bu}  \label{Leibpart5}\\
& &\partial_{\bu} v= r^{2} v \partial_{\bu} +I
\ena
\eqa
& & (P_A)^{ab}_{~~cd} \partial_b \partial_a=0 \label{partialcomm2}\\
& & \partial_b \partial_{\bu} - {q_b \over r^2} \partial_{\bu}
\partial_{b}=0
\ena
 {\bf Conjugation:}
\eq
(x^a)^*=\Dcal{a}{b}x^b,~(dx^a)^*=\Dcal{a}{b}dx^b,~(\partial_a)^*=
-r^Nd^{-1}_a \Dcal{b}{a}  \partial_b
\en
\eq
v^*=v,~ (dv)^*=dv,~(\partial_{\bu})^*=u-\partial_{\bu}
\en
\eq
\overline{r}=r^{-1},~~\overline{q_a}={1 \over q_a} \mbox{for $a
\not=n,n+1$},~~\overline{q_n}={1\over q_{n+1}}
\en
\sk\sk
\noi {\large \bf Table 3: the reduced $SO_{q,r}(N)$ -bicovariant
$x^a,\partial_a,dx^a$
algebra}
\eq
\PA{ab}{cd} x^c x^d=0
\en
\eq
 x \otimes dx=  r\Rha (dx \otimes x)
\en
\eq
dx \we dx =-r {\hat R} (dx \we dx)
\en
\eq
\partial_c x^b= r \Rha^{be}_{~~cd} x^d \partial_e+
\delta^b_c I
\en
\eq
 \PA{ab}{cd} \partial_b \partial_a=0
\en
\sk
{\bf Conjugation:}
\eq
(x^a)^*=\Dcal{a}{b}x^b,~(dx^a)^*=\Dcal{a}{b}dx^b,~(\partial_a)^*=
-r^Nd^{-1}_a \Dcal{b}{a}  \partial_b
\en
\sk\sk

%%%%%%%%%%%%%%%%%%%%%%%%%%%%%%%%%
\app{X, dX, P commutations for the  D=4 quantum Minkowski space }
%%%%%%%%%%%%%%%%%%%%%%%%%%%%%%%%%
\sk
{}~~~~$X$ real and $P$ hermitean: $X^*=X,P^*=P$
\sk
{\bf Parameters }
\sk
Two parameters:  $r$, $q \equiv q_{12}$, with
\eq
|r|=1,~~{q\over r} \in {\bf R}~~\Rightarrow
\overline{r}=r^{-1},~\overline{q}={q\over r^2}
\en
\sk
{\bf Definitions}
\sk
\eq
\la \equiv r-r^{-1},~~\lati \equiv {r\over q}-{q\over r},
{}~~\mu \equiv r+r^{-1},~~\muti \equiv {r\over q}+{q\over r}
\en

{\bf  XX commutations}
\eqa
& & X^2 X^1={1 \over {\mu^2+\lati^2}}[\mu \muti X^1 X^2 + \mu\la X^2
X^4
+i\lati\muti
          X^3 X^4 -i\lati\la X^1 X^3] \nonumber \\
& & X^3 X^1={1 \over {\mu^2+\lati^2}}[\mu \muti X^1 X^3 + \mu\la X^3
X^4
-i\lati\muti
          X^2 X^4 +i\lati\la X^1 X^2] \nonumber \\
& & X^4 X^1 = X^1 X^4 + {\la \over 2} (X^2 X^2 + X^3 X^3) \nonumber
\\
& & X^3 X^2 = X^2 X^3 \nonumber \\
& & X^4 X^2 = {1 \over {\mu^2+\lati^2}}[\mu \muti X^2 X^4-\mu\la X^1
X^2
-i\lati\muti
          X^1 X^3 - i\lati\la X^3 X^4] \nonumber \\
& & X^4 X^3 = {1 \over {\mu^2+\lati^2}}[\mu \muti X^3 X^4-\mu\la X^1
X^3
+i\lati\muti
          X^1 X^2  + i\lati\la X^2 X^4]  \nonumber
\ena
\sk

{\bf X dX commutations}
\eqa
& & X^1 dX^1={1\over 4} (r^{-2}+3r^2) (dX^1)X^1-{\la^2\over
4}[(dX^4)X^1-(dX^1)X^4]
\nonumber\\
& &~~~~~~~~~~~~~       - {\la \over 2} [(dX^2)X^2+(dX^3)X^3] +
{1\over
4}(r^2-r^{-2}) (dX^4)X^4 \nonumber \\
& & X^1 dX^2={1\over 2} r\muti (dX^2)X^1+{r\la \over 2}
d[X^1-X^4]X^2
        -{i\over 2} r\lati (dX^3)X^4 \nonumber \\
& & X^1 dX^3={1\over 2} r\muti (dX^3)X^1+{r\la \over 2}
d[X^1-X^4]X^3
        +{i\over 2} r\lati (dX^2)X^4 \nonumber \\
& & X^1 dX^4=[{1\over 4}(r^2-r^{-2})+1] (dX^4) X^1+{\la^2\over
4}[(dX^1)X^1-(dX^4)X^4]
\nonumber\\
& &~~~~~~~~~~~~~       - {\la \over 2} [(dX^2)X^2+(dX^3)X^3]
+[{1\over
4}(3r^2+r^{-2})-1]
(dX^1)X^4 \nonumber \\
& & X^2 dX^1={1\over 2} r\muti (dX^1)X^2+{r\la \over 2} dX^2
(X^1+X^4)
        +{i\over 2} r\lati (dX^4)X^3 \nonumber \\
& & X^2 dX^2={r\mu \over 2} (dX^2)X^2 -{r\la \over 2}(dX^3)X^3 -{\la
\over 2}
 d[X^1-X^4] (X^1+X^4) \nonumber \\
& & X^2 dX^3 = {r\mu \over 2} (dX^3)X^2 + {r\la \over 2} (dX^2)X^3
\nonumber \\
& & X^2 dX^4={1\over 2} r\muti (dX^4)X^2+{r\la \over 2} dX^2
(X^1+X^4)
        +{i\over 2} r\lati (dX^1)X^3 \nonumber \\
& & X^3 dX^1={1\over 2} r\muti (dX^1)X^3+{r\la \over 2} dX^3
(X^1+X^4)
        -{i\over 2} r\lati (dX^4)X^2 \nonumber \\
& & X^3 dX^2 = {r\mu \over 2} (dX^2)X^3 + {r\la \over 2} (dX^3)X^2
\nonumber \\
& & X^3 dX^3={r\mu \over 2} (dX^3)X^3 -{r\la \over 2}(dX^2)X^2 -{\la
\over 2}
 d[X^1-X^4] (X^1+X^4) \nonumber \\
& & X^3 dX^4={1\over 2} r\muti (dX^4)X^3+{r\la \over 2} dX^3
(X^1+X^4)
        -{i\over 2} r\lati (dX^1)X^2 \nonumber \\
& & X^4 dX^1=[{1\over 4}(r^2-r^{-2})+1] (dX^1) X^4-{\la^2\over
4}[(dX^1)X^1-(dX^4)X^4]
\nonumber\\
& &~~~~~~~~~~~~~       + {\la \over 2} [(dX^2)X^2+(dX^3)X^3]
+[{1\over
4}(3r^2+r^{-2})-1]
(dX^4)X^1 \nonumber \\
& & X^4 dX^2={1\over 2} r\muti (dX^2)X^4-{r\la \over 2}
d[X^1-X^4]X^2
        -{i\over 2} r\lati (dX^3)X^1 \nonumber \\
& & X^4 dX^3={1\over 2} r\muti (dX^3)X^4-{r\la \over 2}
d[X^1-X^4]X^3
        +{i\over 2} r\lati (dX^2)X^1 \nonumber \\
& & X^4 dX^4={1\over 4} (r^{-2}+3r^2) (dX^4)X^4+{\la^2\over
4}[(dX^4)X^1-(dX^1)X^4]
\nonumber\\
& &~~~~~~~~~~~~~       + {\la \over 2} [(dX^2)X^2+(dX^3)X^3] +
{1\over 4}(r^2-r^{-2}) (dX^1)X^1 \nonumber
\ena
\sk

{\bf  dX dX commutations}
\sk
(Products between $dX$ are exterior (wedge) products)
\eqa
& & dX^1  dX^1 = 0  \nonumber \\
& & dX^1 dX^2= {1 \over \mu^2-\lati^2}[-\mu\muti dX^2  dX^1 + \mu\la
dX^4 dX^2
-i\lati\muti dX^4 dX^3 -i\la\lati dX^3 dX^1]  \nonumber \\
& & dX^1 dX^3 = {1 \over \mu^2-\lati^2}[-\mu\muti dX^3  dX^1 +
\mu\la
dX^4 dX^3
+i\lati\muti dX^4 dX^2 +i\la\lati dX^2 dX^1]  \nonumber \\
& & dX^1 dX^4=-dX^4 dX^1  \nonumber \\
& & dX^2 dX^2= -{\la\over 2} dX^4 dX^1  \nonumber \\
& & dX^2 dX^3=-dX^3 dX^2  \nonumber \\
& & dX^2 dX^4= {1 \over \mu^2-\lati^2}[-\mu\muti dX^4  dX^2 - \mu\la
dX^2 dX^1
+i\lati\muti dX^3 dX^1 -i\la\lati dX^4 dX^3]  \nonumber \\
& & dX^3 dX^3=  -{\la\over 2} dX^4 dX^1  \nonumber \\
& & dX^3 dX^4= {1 \over \mu^2-\lati^2}[-\mu\muti dX^4  dX^3 - \mu\la
dX^3 dX^1
-i\lati\muti dX^2 dX^1 +i\la\lati dX^4 dX^2]  \nonumber \\
& & dX^4 dX^4 = 0  \nonumber
\ena

\sk

{\bf  P X commutations }
\sk
Defining
\eqa
{}~~~~A&\equiv& -{r\la \over 4} [ (X^1 - X^4)(P_1 -P_4)
+ r^2 (X^1 + X^4)(P_1 + P_4) + 2r(X^2P_2+X^3P_3)]\nonumber\\
{}~~~~B&\equiv&  -{r\la \over 4} [ (X^1 - X^4)(P_1 -P_4)
 - r^2 (X^1 + X^4)(P_1 + P_4) - 2r(X^2P_2+X^3P_3)]\nonumber
\ena
\indent  the commutations are:
\eqa
& & P_1X^1-X^1P_1 +A= -{1\over 2} i \hbar r^2 \mu
\nonumber \\
& &P_1X^2-{r\muti\over 2} X^2P_1- {r \over 2} [- \la (X^1+X^4)P_2
+ i
\lati
X^3P_4]=0
 \nonumber \\
& &P_1X^3-{r\muti\over 2} X^3P_1- {r \over 2} [- \la (X^1+X^4)P_3- i
\lati
X^2P_4]=0
 \nonumber \\
& & P_1X^4- r^2 X^4P_1 +B -{r\la} X^1P_4 = {1\over 2} i\hbar r^2 \la
\nonumber \\
& &P_2X^1-{r\muti\over 2} X^1P_2- {r \over 2} [ - \la X^2 (P_1+P_4)+
i \lati
X^4P_3]=0
 \nonumber \\
& &P_2X^2-{r\mu\over 2} X^2P_2- {r \la \over 2} [ X^3P_3 + r
(X^1+X^4)
(P_1+P_4)]=-i\hbar r^2
 \nonumber \\
& &P_2X^3-{r\mu\over 2} X^3P_2 + {r \la \over 2} X^2P_3 = 0
 \nonumber \\
& &P_2X^4-{r\muti\over 2} X^4 P_2- {r \over 2} [  \la X^2 (P_1+P_4)+
i \lati
X^1P_3]=0
 \nonumber \\
& &P_3X^1-{r\muti\over 2} X^1P_3- {r \over 2} [ - \la X^3 (P_1+P_4) -
i \lati
X^4P_2]=0
 \nonumber \\
& &P_3X^2-{r\mu\over 2} X^2P_3 + {r \la \over 2} X^3P_2 = 0
 \nonumber \\
& &P_3X^3-{r\mu\over 2} X^3P_3- {r \la \over 2} [ X^2P_2 + r
(X^1+X^4)
(P_1+P_4)]=-i\hbar r^2
 \nonumber \\
& &P_3X^4-{r\muti\over 2} X^4P_3- {r \over 2} [  \la X^3 (P_1+P_4) -
i \lati
X^1P_2]=0
 \nonumber \\
& & P_4X^1- r^2 X^1P_4  +B -{r\la} X^4P_1 = {1\over
2} i
\hbar r^2 \la
\nonumber \\
& &P_4X^2-{r\muti\over 2} X^2P_4- {r \over 2} [ \la (X^1+X^4)P_2+ i
\lati
X^3P_1]=0
 \nonumber \\
& &P_4X^3-{r\muti\over 2} X^3P_4- {r \over 2} [ \la (X^1+X^4)P_3 - i
\lati
X^2P_1]=0
 \nonumber \\
& & P_4X^4-X^4P_4 +A = -{1\over 2} i \hbar r^2 \mu
\nonumber
\ena
\sk

{\bf P P commutations}

\eqa
& & P^2 P^1={1 \over {\mu^2+\lati^2}} [\mu \muti P^1 P^2- \mu\la
P^2
P_4
- i \lati \muti P_3 P_4 - i \lati \la P_1 P_3] \nonumber \\
& & P^3 P^1={1 \over {\mu^2+\lati^2}} [\mu \muti P^1 P^3- \mu\la P^3
P_4
+ i \lati \muti P_2 P_4 + i \lati \la P_1 P_2] \nonumber \\
& & P^4 P^1 = P^1 P^4 - {\la \over 2} (P_2P_2+P_3P_3) \nonumber\\
& & P^3 P^2 = P^2 P^3 \nonumber \\
& & P^4 P^2={1 \over {\mu^2+\lati^2}} [\mu \muti P^2 P^4 + \mu\la
P^1P_2
+ i \lati \muti P_1 P_3 - i \lati \la P_3 P_4] \nonumber \\
& & P^4 P^3={1 \over {\mu^2+\lati^2}} [\mu \muti P^3 P^4 + \mu\la
P^1P_3
- i \lati \muti P_1 P_2 + i \lati \la P_2 P_4] \nonumber
\ena

%%%%%%%%%%%%%%%%%%%%%%%%%%%%%%%%%
\sect{$R$ matrix of orthogonal  $q$-groups: properties}

Let $\Rh$ be  the matrix defined by
$\Rhat{ab}{cd} \equiv \R{ba}{cd}$.
\sk
Characteristic equation and  projector decomposition:
\eq
(\Rh-rI)(\Rh+r^{-1}I)(\Rh- r^{1-N} I)=0 \label{cubic}
\en
\eq
\Rh-\Rh^{-1}=(r-r^{-1}) (I-K) \label{extrarelation}
\en
\eq
\Rh=r P_S - r^{-1} P_A+ r^{1-N}P_0  \label{RprojBCD}
\en
with
\eq
\begin{array}{ll}
&P_S={1 \over {r+\rminus}} [\Rh+\rminus I-(\rminus+
r^{1-N})P_0]\\
&P_A={1 \over {r+\rminus}} [-\Rh+rI-(r- r^{1-N})P_0]\\
&P_0= Q_N(r) K\\
&Q_N(r) \equiv (C_{ab} C^{ab})^{-1}={{1-r^{-2}} \over
{(1- r^{-N})(1+ r^{N-2})}}~,~~~~
K^{ab}_{~~cd}=C^{ab} C_{cd}\\
&I=P_S+P_A+P_0
\end{array}
\label{projBCD}
\en
Other properties involving the $q$-metric:
\eq
C_{ab} \Rhat{bc}{de} = \Rhatinv{cf}{ad} C_{fe} ,~~~
 \Rhat{bc}{de} C^{ea}=C^{bf} \Rhatinv{ca}{fd} \label{crc}
\en
\eq
C_{ab}\Rhat{ab}{cd}= r^{1-N}C_{cd} ,~~~
C^{cd}\Rhat{ab}{cd}= r^{1-N}C^{ab}  \label{CR}
\en
\noi The identities (\ref{crc})
 hold also for $\Rh \rightarrow \Rh^{-1}$.
\sk
{\bf Acknowledgments}
\sk
This work is  supported in part by EEC under TMR contract
FMRX-CT96-0045;  by the Director,
Office of Energy Research, Office af High Energy  and Nuclear
Physics,
Division of High Energy Physics
of the U.S. Department of Energy under Contract DE-AC03-76SF00098
and by the National Science Foundation under grant PHY-95-14797.
P.A. is supported by  a research fellowship of
Fondazione Angelo Della Riccia.

%%%%%%%%%%%%%%%%%%%%%%%%%%%%%%%%%

\vfill\eject
\end{document}